\documentclass[times, authoryear, 12pt, twoside]{elsarticle}
\usepackage{lmodern}
\usepackage{geometry}
\usepackage[labelfont=bf]{caption}
\usepackage{amsmath,amsfonts,amssymb,amsthm,color,epsfig,graphicx,url, eurosym, latexsym}
\usepackage{appendix}
\usepackage[utf8]{inputenc}
\usepackage[english]{babel}

\urldef\myurl\url{www.ualberta.ca/~dwiens/}
\newtheorem{theorem}{Theorem}

\newtheorem{lemma}{Lemma}
\newtheorem{remark}{Remark}

\makeatletter
\def\ps@pprintTitle{  \let\@oddhead\@empty
  \let\@evenhead\@empty
  \def\@oddfoot{\reset@font\hfil\thepage\hfil}
  \let\@evenfoot\@oddfoot
}

\begin{document}

\def\T{{ \mathrm{\scriptscriptstyle T} }}

\date{\today }

\input{amssym.def}

\begin{frontmatter}
\title{Minimax Robust Designs for M-Estimated Models}

\author[A1]{Rui Hu\corref{mycorrespondingauthor1}}
\address[A1]{Mathematics \& Statistics,
	MacEwan University,
	Edmonton, Canada,  T5J 4S2}

\author[A2]{Douglas P. Wiens\corref{mycorrespondingauthor2}}
\address[A2]{Mathematical \& Statistical Sciences,
	University of Alberta,
	Edmonton, Canada,  T6G 2G1}

\cortext[mycorrespondingauthor1]{E-mail: \url{hur3@macewan.ca}. } 
\cortext[mycorrespondingauthor2]{E-mail: \url{doug.wiens@ualberta.ca}. Computing code is at {\myurl}.} 

\begin{abstract}
Experimental designs that are minimax in the presence of
model misspecifications have been constructed so as to minimize the maximum,
over classes of alternate response models, of the integrated mean squared
error of the predicted values. The theory to date has focussed almost exclusively on Least Squares estimates. Here we extend this theory to designs tailored for M-estimation of parameters, thus obtaining additional robustness against outlying, or otherwise erroneous, responses. We show that, subject to a minor change in a tuning constant, designs optimal for Least Squares remain so asymptotically for M-estimation, under squared error loss. We argue that even this minor change should be ignored, and the tuning constant chosen in an \textit{ad hoc} but sensible manner that does not depend on which M-estimate is being employed. A rather surprising additional result is that our designs and estimates, derived under an assumption of i.i.d.\ errors, are also robust, in a minimax sense, against broad classes of correlation structures. We also investigate designs under a more robust loss function, less influenced by outlying responses than squared error.
\end{abstract}

\begin{keyword} 
asymptotics \sep
finite design space \sep
Integrated M-loss \sep
misspecified model \sep
regression design \sep
robustness against dependence. 
\end{keyword}
\end{frontmatter}

\section{Introduction and summary \label{section: intro}}

The theory and practice of robustness of design, for possibly misspecified
response functions, is well-developed as it applies to cases in which
parameter estimation is to be carried out by Least Squares (\textsc{ls}).
See \cite{w15} for a review. An investigator seeking model robustness might
naturally be concerned as well with robustness against outlying data points,
or more generally against a misspecified data-generating probability
distribution, and hence seek M-estimates of the parameters. Indeed, as
pointed out by \cite{m97}, `a robust inference will profit from an optimal
design and ... an optimal design is more reasonable if it allows also a
robust analysis of the data'. See also \cite{mm01} for the special case of $%
L_1$ estimation.

\cite{w94,w96a} studied this design problem for quite limited classes of
approximate responses on continuous design spaces and obtained asymptotic
results under rather restrictive conditions. \cite{ww10} carried out a small
simulation study and found that there was little apparent dependence of the
designs on the method of estimation.

In this article we strengthen and extend these findings. We show that, if
the parameters of the assumed model are to be estimated by ordinary
M-estimation, then designs optimal, in the sense of minimizing the maximum
integrated mean squared error (\textsc{imse}) of the predictions, for \textsc{ls} remain so,
asymptotically, for M-estimation.

We note that M-estimates offer protection against outlying responses, but
not necessarily against outliers in the factor space. This latter type of
protection is furnished by Generalized M-estimation. But since our focus
centres on design points chosen by the experimenter, such outlyingness is
not an issue.

In \S \ref{section: asymptotics} of this article we present the asymptotic
theory on which our design problem will be based. The challenges there are
somewhat unique, since we do not assume that the fitted response model is
the correct one, and must allow for a broad class of alternatives. The
ensuing Theorem \ref{thm:main} seems to be the most explicit and
comprehensive response to this challenge, in the literature. In \S \ref{section: minimax theory} we address the design problem. We show
that the optimally robust designs, under squared error loss, depend on the anticipated method of
estimation only through a tuning constant. We later argue that although this tuning constant depends on the 
anticipated method of estimation through unknown parameters,
this dependence is so slight that it should be ignored and the tuning
constant chosen in an \textit{ad hoc} but sensible manner -- as the relative
emphasis placed by the designer on bias reduction versus variance reduction.
The result is then that the designs are completely independent of which
M-estimate is to be used. 

In \S \ref{section: dependence} we present a surprising (to us) result under which the designs
and estimates presented here are also minimax robust against broad classes of correlation
structures; a specific case is that in which the random errors are
equicorrelated.

In \S \ref{section: IMloss} we consider
a more robust but computationally much more expensive alternative to \textsc{imse} -- `Integrated M-loss' (\textsc{im}), tailored to the robust estimation method. Whereas the `max' part of the minimax problem addressed in \S \ref{section: minimax} is carried out analytically, and in \S \ref{section: construction} we present a very efficient algorithm for the minimization, the construction of minimax \textsc{im} designs requires both the `max' and `min' components of the problem to be carried out numerically. Algorithms for the construction of  minimax \textsc{im} designs are discussed in 
\S \ref{section: IMloss}, with examples of both types of designs in in \S \ref{section: examples}. 

Proofs are in the Appendix. The \textsc{matlab} code used to prepare the
examples is available on the second author's personal website.

\section{Minimax robustness of design \label{section: minimax}}

\subsection{Asymptotic theory \label{section: asymptotics}}

Our minimax design problem is phrased in terms of an approximate regression
response 
\begin{equation}
E\left[ Y\left( \boldsymbol{x}\right) \right] \approx \boldsymbol{f}^{\T }\left( \boldsymbol{x}\right) \boldsymbol{%
\theta },  \label{approx}
\end{equation}%
for $p$ regressors $\boldsymbol{f}$, each functions of vectors $\boldsymbol{x%
}$ of $q$ independent variables, ranging over a finite \textit{design space} 
$\chi=\left\{ \boldsymbol{x}_{1},...,\boldsymbol{x}_{N}\right\} \subset 
\mathbb{R}^{q}$, and for a parameter vector $\boldsymbol{\theta }_{p\times 1}
$.

Since (\ref{approx}) is an approximation, the interpretation of $\boldsymbol{%
\theta }$ is unclear; we \textit{define }this target parameter by 
\begin{equation}
\boldsymbol{\theta }_{0}=\arg \min_{\boldsymbol{\eta }}\sum_{i=1}^{N}\left( E%
\left[ Y\left( \boldsymbol{x}_{i}\right) \right] -\boldsymbol{f}^{\T }\left( \boldsymbol{x}_{i}\right) \boldsymbol{\eta }%
\right) ^{2}.  \label{theta def}
\end{equation}%
Equivalently, and with model error $\tau \left( \boldsymbol{x}\right) 
\overset{def}{=}E\left[ Y\left( \boldsymbol{x}\right) \right] -\boldsymbol{f}^{\T}\left( \boldsymbol{x}\right) 
\boldsymbol{\theta }_{0}$, we have 
\begin{equation}
\sum_{i=1}^{N}\boldsymbol{f}\left( \boldsymbol{x}_{i}\right) \tau \left( 
\boldsymbol{x}_{i}\right) =\boldsymbol{0}_p.   \label{orthogonality}
\end{equation}%
Assuming that $\mathcal{X}$ is rich enough that the matrix\ $\boldsymbol{A}%
\overset{def}{=}\sum_{i=1}^{N}\boldsymbol{f}\left( \boldsymbol{x}_{i}\right) 
\boldsymbol{f}^{\T}\left( \boldsymbol{x}%
_{i}\right) $ is invertible, the parameter defined by (\ref{theta def}) and (%
\ref{orthogonality}) is unique. In order that variance and bias remain of a
comparable order asymptotically, we bound the approximation error in (\ref%
{approx}) by assuming that%
\begin{equation}
\sum_{i=1}^{N}\tau ^{2}\left( \boldsymbol{x}_{i}\right) \leq \kappa ^{2}/n,
\label{bound}
\end{equation}%
for a constant $\kappa $. Here $n$ is the total number of experimental runs
to be carried out.

Our model is thus given by $E\left[ Y\left( \boldsymbol{x}\right) \right] =%
\boldsymbol{f}^{\T}\left( \boldsymbol{x}%
\right) \boldsymbol{\theta }_{0}+\tau \left( \boldsymbol{x}\right) $, for an
unknown model error $\tau \left( \boldsymbol{\cdot }\right) $ constrained by
(\ref{orthogonality}) and (\ref{bound}). Let $\Upsilon $ be the class of
such functions $\tau \left( \boldsymbol{\cdot }\right) $.

At such values of $\boldsymbol{x}$, observations $Y\left( \boldsymbol{x}%
\right) $ may be made, with additive random error: $Y\left( \boldsymbol{x}%
\right) =E\left[ Y\left( \boldsymbol{x}\right) \right] + \varepsilon ,$ for
i.i.d., symmetrically distributed errors $\varepsilon $ with finite
variance. We make no further distributional assumptions. In particular, the
distribution of $\varepsilon $ may be different from that assumed by the
experimenter, thus calling for the type of robustness furnished by
M-estimation.

Thus, given observations $\left\{ Y_{j}\left( \boldsymbol{x}_{i}\right)
|j=1,...,n_{i}\right\} $, with $n_{i}$ replicates at distinct points $%
\boldsymbol{x}_{i}$, we suppose that $\boldsymbol{\theta }_{0}$ will be
estimated by M-estimation with an auxiliary estimate of scale. For a
function $\psi $ with properties as in C1) below, the estimate satisfies 
\begin{equation}
\sum_{i,j}\psi \left( \frac{Y_{j}\left( \boldsymbol{x}_{i}\right) -%
\boldsymbol{f}^{\T}\left( \boldsymbol{x}%
_{i}\right) \boldsymbol{\hat{\theta}}_{n}}{\hat{\sigma}_{n}}\right) 
\boldsymbol{f}\left( \boldsymbol{x}_{i}\right)=\boldsymbol{0}_p ,
\label{regression}
\end{equation}%
where $\hat{\sigma}_{n}$ is a consistent estimate of scale, for instance 
\begin{equation}
\hat{\sigma}_{n}=\left. \text{\textsc{median}}\left\{ \left\vert Y_{j}\left( 
\boldsymbol{x}_{i}\right) -\boldsymbol{f}^{{\T} }\left( \boldsymbol{x}_{i}\right) \boldsymbol{\hat{\theta}}_{n}\right\vert
\right\} \right/ \Phi ^{-1}\left( .75\right),  \label{scale}
\end{equation}%
which is consistent for the standard deviation $\sigma $ if the data are
Normal. Under the mild conditions of \cite{rs86}, $\boldsymbol{\hat{\theta}}%
_{n}$ has the same asymptotic properties as if $\hat{\sigma}_{n}$ were
replaced by $\sigma $. Thus, with $\psi _{\sigma }\left( x\right) \overset{%
def}{=}\psi \left( x/\sigma \right) $, we have that 
\begin{equation}
\sum_{i,j}\psi _{\sigma }\left( Y_{j}\left( \boldsymbol{x}_{i}\right) -%
\boldsymbol{f}^{\T}\left( \boldsymbol{x}%
_{i}\right) \boldsymbol{\hat{\theta}}_{n}\right) \boldsymbol{f}\left( 
\boldsymbol{x}_{i}\right)=\boldsymbol{0}_p,  \label{m-eqn}
\end{equation}%
and we replace $\sigma $ by $\hat{\sigma}_{n}$ in the applications. For the
computations one would iterate between (\ref{regression}) and (\ref{scale}).
See \cite{h81} for further details of M-estimation, \cite{hw17} for
applications to testing and \cite{h18} for applications in spatial sampling.

There is a rich literature concerning minimax variance M-estimation when the
true error distribution is assumed to vary over a neighbourhood of the
assumed distribution -- see for instance \cite{h64}, \cite{w86}, \cite{cw85}%
, \cite{cw89} for examples of solutions in various such neighbourhoods. We
do not pursue this further, and leave the choice of a specific M-estimator
up to the experimenter.

By a \textit{design} we mean a probability mass function $\xi$ on $\mathcal{X%
}$, to be approximated as necessary by an implementable, $n$-point design $%
\xi_{n}$, with $\xi_{n}(\boldsymbol{x}_{i})=n_{i}/n$. For any design $\xi$
define 
\begin{equation*}
\boldsymbol{M}_{0}\left( \xi \right) =\sum_{i=1}^{N}\boldsymbol{f}\left( 
\boldsymbol{x}_{i}\right) \boldsymbol{f}^{\T}\left( \boldsymbol{x}_{i}\right) \xi \left( \boldsymbol{x}_{i}\right) \text{
and }\boldsymbol{b}_{0}\left( \xi \right) =\sum_{i=1}^{N}\boldsymbol{f}%
\left( \boldsymbol{x}_{i}\right) \tau\left( \boldsymbol{x}_{i}\right) \xi
\left( \boldsymbol{x}_{i}\right) .
\end{equation*}

We define $\sigma _{M}^{2}=\left. E\left[ \psi _{\sigma }^{2}\left(
\varepsilon \right) \right] \right/ \left( E\left[ \psi _{\sigma }^{\prime
}\left( \varepsilon \right) \right] \right) ^{2}$ (assumed finite), and for
a design $\xi _{n}$ we set $\boldsymbol{M}_{n}=\boldsymbol{M}_{0}\left( \xi
_{n}\right) $, and $\boldsymbol{b}_{n}=\boldsymbol{b}_{0}\left( \xi
_{n}\right) $. In the Appendix we prove Theorem \ref{thm:main}, stated at
the end of this section, which makes precise the asymptotic normality of $%
\boldsymbol{\hat{\theta}}_{n}-\boldsymbol{\theta }_{0}$:%
\begin{equation}
\boldsymbol{\hat{\theta}}_{n}-\boldsymbol{\theta }_{0}\sim AN\left( 
\boldsymbol{M}_{n}^{-1}\boldsymbol{b}_{n},\left( \sigma _{M}^{2}/n\right) 
\boldsymbol{M}_{n}^{-1}\right) .  \label{AN}
\end{equation}%
\noindent For this we define $\boldsymbol{\tilde{\theta}}_{n}=\boldsymbol{%
\theta }_{0}+\boldsymbol{M}_{n}^{-1}\boldsymbol{b}_{n}$, and will show that $%
\sqrt{n}\left( \boldsymbol{\hat{\theta}}_{n}-\boldsymbol{\theta }_{0}-%
\boldsymbol{M}_{n}^{-1}\boldsymbol{b}_{n}\right) =\sqrt{n}\left( \boldsymbol{%
\hat{\theta}}_{n}-\boldsymbol{\tilde{\theta}}_{n}\right) $ is asymptotically
normal with mean $\boldsymbol{0}_p$.

For the proof of Theorem \ref{thm:main} we make the following assumptions
and definitions. We denote the smallest and largest eigenvalues of a matrix
by $ch_{\min }$ and $ch_{\max }$ respectively.

\begin{description}
\item[C1)] The function $\psi \left( \cdot \right) $ is weakly increasing,
twice differentiable and odd. We define $m_{1}=\sup_t \psi _{\sigma
}^{\prime }\left(t\right) $ and $m_{2}=\sup_t|\psi _{\sigma }^{\prime \prime
}\left(t\right) |$, both assumed finite.

\item[C2)] The sequence $A_{n}=ch_{\min }\boldsymbol{M}_{n}$ is bounded away
from zero.

\item[C3)] The sequence of designs has a weak limit $\xi _{\ast }$: $\xi _{n}%
\overset{d}{\rightarrow }\xi _{\ast }$ as $n\rightarrow \infty $.

\begin{remark}
\label{remark: assumptions} Assumption C1) is standard, but precludes the
use of redescenders (\cite{m04}). It implies that $E\left[ \psi \left(
\varepsilon \right) \right] =E\left[ \psi ^{\prime \prime }\left(
\varepsilon \right) \right] =0$, since the distribution of $\varepsilon $ is
symmetric, and that $E\left[ \psi ^{\prime }\left( \varepsilon \right) %
\right] >0$, else $\psi \left( \cdot \right) $ is identically zero. For C2),
that $ch_{\min }\boldsymbol{M}_{n} >>0$ ensures that $\boldsymbol{M}_{n}$
remains invertible, and is a natural property of the design. For C3), that
the designs have a weak limit asserts only that the design weights $\xi_{n}(%
\boldsymbol{x}_{i})$ converge, which is a requirement of our numerical
algorithm, discussed in \S \ref{section: construction}. The limit will be
used in the proof of Theorem \ref{thm:main} when we invoke Theorem 5.9 of 
\cite{vw98}, upon which we base our proof that $\boldsymbol{\hat{\theta}}%
_{n}-\boldsymbol{\tilde{\theta}}_{n}$ is a $\sqrt{n}$-consistent estimate of
zero.
\end{remark}
\end{description}

\begin{theorem}
\label{thm:main}With notation as above, we have that, as $n\rightarrow
\infty $, 
\begin{equation*}
\sqrt{n}\boldsymbol{M}_{n}^{1/2}\left( \boldsymbol{\hat{\theta}}_{n}-%
\boldsymbol{\theta }_{0}-\boldsymbol{M}_{n}^{-1}\boldsymbol{b}_{n}\right) 
\overset{d}{\longrightarrow }N\left( \boldsymbol{0}_p,\sigma _{M}^{2}%
\boldsymbol{I}_{p}\right) .
\end{equation*}
\end{theorem}

\subsection{Minimax design theory \label{section: minimax theory}}

We define our loss in terms of the Integrated Mean Squared Error of the
predictors $\hat{Y}\left( \boldsymbol{x}\right) =\boldsymbol{f}^{\T}\left( \boldsymbol{x}\right) \boldsymbol{\hat{\theta%
}}_{n}$: 
\begin{equation}
\text{\textsc{imse}}\left( \xi |\tau \right) =\sum_{\boldsymbol{x}\mathbb{\
\in }\chi}E\left\{ \left( E\left[ Y\left( \boldsymbol{x}\right) \right] -%
\hat{Y}\left( \boldsymbol{x}\right) \right) ^{2}\right\} =\sum_{\boldsymbol{x%
}\mathbb{\in }\chi}E\left\{ \left( \tau \left( \boldsymbol{x}\right) -%
\boldsymbol{f}^{\T}\left( \boldsymbol{x}%
\right) \left( \boldsymbol{\hat{\theta}}_{n}\left( \xi |\tau \right) -%
\boldsymbol{\theta }_{0}(\tau)\right) \right) ^{2}\right\} .  \label{IM}
\end{equation}%
We will maximize an asymptotic approximation to \textsc{imse}$\left( \xi
|\tau \right) $ over $\tau $ satisfying (\ref{orthogonality}) and (\ref%
{bound}), and then find designs minimizing this maximum. The approximation
uses the asymptotic mean and covariance given in (\ref{AN}). See \cite{l99}
for a discussion  limiting versus asymptotic moments.

Expanding (\ref{IM}) and using (\ref{orthogonality}) results in 
\begin{equation}
\text{\textsc{imse}}\left( \xi|\tau \right)=\sum_{i=1}^{N}\tau ^{2}\left( 
\boldsymbol{x}_{i}\right) +\sum_{i=1}^{N}\boldsymbol{f}^{\T}\left( \boldsymbol{x}_{i}\right) E\left[ \left( 
\boldsymbol{\hat{\theta}}_{n}(\xi|\tau)-\boldsymbol{\theta }%
_{0}(\tau)\right)\left( \boldsymbol{\hat{\theta}}_{n}\left( \xi |\tau
\right) -\boldsymbol{\theta }_{0}(\tau)\right) ^{\T}\right] \boldsymbol{f}\left( \boldsymbol{x}%
_{i}\right).  \notag
\end{equation}
Substituting the asymptotic moments gives the asymptotic approximation 
\begin{eqnarray}
\text{\textsc{aimse}}\left( \xi_n |\tau \right)&=&\sum_{i=1}^{N}\tau
^{2}\left( \boldsymbol{x}_{i}\right) +\sum_{i=1}^{N}\boldsymbol{f}^{\T}\left( \boldsymbol{x}_{i}\right) \left[ 
\boldsymbol{M}_{n}^{-1}\boldsymbol{b}_{n}\boldsymbol{b}_{n}^{\T}\boldsymbol{M}_{n}^{-1}+\left( \sigma
_{M}^{2}/n\right) \boldsymbol{M}_{n}^{-1}\right] \boldsymbol{f}\left( 
\boldsymbol{x}_{i}\right)  \notag \\
&=&\sum_{i=1}^{N}\tau ^{2}\left( \boldsymbol{x}_{i}\right) +\boldsymbol{b}%
_{n}^{\T}\boldsymbol{M}_{n}^{-1}\mathbf{A}%
\boldsymbol{M}_{n}^{-1}\boldsymbol{b}_{n}+\left( \sigma _{M}^{2}/n\right) tr%
\boldsymbol{AM}_{n}^{-1} \text{.}  \label{IMSEn}
\end{eqnarray}

For the maximization it is convenient to introduce an orthogonal basis for
the space of regressors. Define $\boldsymbol{F}_{N\times p}=\left( 
\boldsymbol{f}\left( \boldsymbol{x}_{1}\right) ,\cdot \cdot \cdot ,
\boldsymbol{f}\left( \boldsymbol{x}_{N}\right) \right)^{\T} $,
and for a design $\xi $ on $\chi $, 
$\boldsymbol{D}\left( \xi \right) =diag\left( \xi _{1},\cdot \cdot \cdot ,\xi _{N}\right) 
$. By the Gram-Schmidt process we can construct a matrix $\boldsymbol{Q}%
_{N\times p}$ whose orthonormal columns form a basis for the column space of 
$\boldsymbol{F}$ -- assumed to be of dimension $p$.

In the Appendix we prove the following theorem.

\begin{theorem}
\label{thm: max loss} For any design for which $\boldsymbol{M}_{0}\left( \xi \right)$ in non-singular, define $p\times p$ matrices 
\begin{equation*}
\boldsymbol{R}\left( \xi \right) =\boldsymbol{Q}^{\T }\boldsymbol{D}%
\left( \xi \right) \boldsymbol{Q},\text{ }\boldsymbol{S}\left( \xi \right) =%
\boldsymbol{Q}^{\T }\boldsymbol{D}^{2}\left( \xi \right) \boldsymbol{Q},%
\text{ }\boldsymbol{U}\left( \xi \right) =\boldsymbol{R}^{-1}\left( \xi
\right) \boldsymbol{S}\left( \xi \right) \boldsymbol{R}^{-1}\left( \xi
\right) .
\end{equation*}%
Then $\max_{\tau \in \Upsilon } \textsc{aimse} \left( \xi _{n}|\tau \right) $
is given by $n^{-1}$ times 
\begin{equation}
J\left( \xi _{n}\right) =\sigma _{M}^{2}tr\boldsymbol{R}^{-1}\left( \xi
_{n}\right) +\kappa ^{2}ch_{\max }\boldsymbol{U}\left( \xi _{n}\right) .
\label{I3}
\end{equation}
\end{theorem}

A further maximization of the \textsc{aimse} is discussed in the next
section, after which, in \S \ref{section: construction}, we minimize this
maximum, thus obtaining the minimax designs.

\subsection{Robustness against dependence \label{section: dependence}}

It is brought out in the proof of Theorem \ref{thm:main} -- see (\ref%
{represent}) of the Appendix -- that the estimate is representable as 
\begin{equation*}
\boldsymbol{\Phi }_{n}\overset{def}{=}\sqrt{n}\boldsymbol{M}_{n}^{1/2}\left( 
\boldsymbol{\hat{\theta}}_{n}-\boldsymbol{\theta }_{0}-\boldsymbol{M}%
_{n}^{-1}\boldsymbol{b}_{n}\right) =\boldsymbol{K}_{p\times n}\boldsymbol{u}+%
\boldsymbol{k}_{p\times 1}+o_{p}\left( 1\right) ,
\end{equation*}%
where $\boldsymbol{K}$ and $\boldsymbol{k}$ are non-random and $\boldsymbol{u%
}=\left( U_{1},...,U_{n}\right) ^{\T} $ with 
$U_{i}=\psi _{\sigma }\left( \varepsilon _{i}\right) \left/ E\left[ \psi
_{\sigma }^{\prime }\left( \varepsilon \right) \right] \right. $. We have up
to now assumed that the errors $\varepsilon _{i}$ are i.i.d., so that as
well the $U_{i}$ are i.i.d., with \textsc{cov}$\left[ \boldsymbol{u}\right]
=\sigma _{M}^{2}\boldsymbol{I}_{n}$ and, asymptotically, \textsc{cov}$\left[ 
\boldsymbol{\Phi }_{n}\right] =\sigma _{M}^{2}\boldsymbol{KK}^{\T}$. We now investigate the effect on the estimate if
the $U_{i}$ are correlated, or are heteroscedastic.

Suppose then that \textsc{cov}$\left[ \boldsymbol{u}\right] =\boldsymbol{C}%
_{n\times n}$, so that the asymptotic covariance structure of $\boldsymbol{%
\Phi }_{n}$ becomes \textsc{cov}$\left[ \boldsymbol{\Phi }_{n}|\boldsymbol{C}%
\right] =\boldsymbol{KCK}^{\T}$. Note that 
\textsc{cov}$\left[ \boldsymbol{\Phi }_{n}|\boldsymbol{C}\right] $ is
non-decreasing in the Loewner ordering: $\boldsymbol{C}_{1}\preceq 
\boldsymbol{C}_{2}\Rightarrow $ \textsc{cov}$\left[ \boldsymbol{\Phi }_{n}|%
\boldsymbol{C}_{1}\right] \preceq $ \textsc{cov}$\left[ \boldsymbol{\Phi }%
_{n}|\boldsymbol{C}_{2}\right] $. If loss is measured by $\mathcal{L}\left( 
\boldsymbol{C}\right) =\phi \left( \text{\textsc{cov}}\left[ \boldsymbol{%
\Phi }_{n}|\boldsymbol{C}\right] \right) $, where $\phi \left( \cdot \right) 
$ is a function, such as the trace, determinant, maximum eigenvalue etc.
that is itself non-decreasing in the Loewner ordering, then $\mathcal{L}%
\left( \boldsymbol{C}\right) $ is non-decreasing in this ordering. In the
particular problem at hand, $\phi \left( \Sigma \right) =\sum_{i=1}^{N}%
\boldsymbol{f}^{\T}\left( \boldsymbol{x}%
_{i}\right) \Sigma \boldsymbol{f}\left( \boldsymbol{x}_{i}\right) $ is the
integrated variance of the predictors when evaluated at $\boldsymbol{\Sigma }%
=$ \textsc{cov}$\left[ \Phi _{n}| \boldsymbol{C}\right] $.

Suppose now that $\left\Vert \cdot \right\Vert _{M}$ is a matrix norm,
induced by a vector norm $\left\Vert \cdot \right\Vert _{V}$, i.e. $%
\left\Vert \boldsymbol{C}\right\Vert _{M}=\sup_{\left\Vert x\right\Vert _{V}%
\text{ }=1}\left\Vert \boldsymbol{Cx}\right\Vert _{V}.$ Special cases are
the spectral radius $\left\Vert \boldsymbol{C}\right\Vert _{E}$ -- this is
the maximum eigenvalue since $\boldsymbol{C}$ is a covariance matrix -- and
the maximum absolute row sum $\left\Vert \boldsymbol{C}\right\Vert
_{1}=\max_{i}\sum_{j}\left\vert c_{ij}\right\vert $.

The following lemma is given in \cite{w25}; we repeat it here for
convenience.

\begin{lemma}
\label{lemma: dependence}For $\eta ^{2}>0$, covariance matrix $\boldsymbol{C}
$ and induced matrix norm $\left\Vert \boldsymbol{C}\right\Vert _{M}$,
define 
\begin{equation*}
\mathcal{C}_{M}=\left\{ \boldsymbol{C}\left\vert {}\right. \boldsymbol{C}%
\succeq 0\text{ and }\left\Vert \boldsymbol{C}\right\Vert _{M}\leq \eta
^{2}\right\} .
\end{equation*}%
For the norm $\left\Vert \mathbf{\cdot }\right\Vert _{E}$ an equivalent
definition is $\mathcal{C}_{E}=\{\boldsymbol{C}\mid 0\preceq \boldsymbol{C}%
\preceq \eta ^{2}\boldsymbol{I}_{n}\}$. Then \noindent (i) in any such class 
$\mathcal{C}_{M}$, $\max_{\mathcal{C}_{M}}\mathcal{L}\left( \boldsymbol{C}%
\right) =\mathcal{L}\left( \eta ^{2}\boldsymbol{I}_{n}\right) $, and
\noindent (ii) if $\mathcal{C}^{\prime }$ $\mathcal{\subseteq C}_{M}$ and $%
\eta ^{2}\boldsymbol{I}_{n}\in \mathcal{C}^{\prime }$, then $\max_{\mathcal{C%
}^{\prime }}\mathcal{L}\left( \boldsymbol{C}\right) =\mathcal{L}\left( \eta
^{2}\boldsymbol{I}_{n}\right) $.
\end{lemma}

A consequence of (i) of this lemma is that if one is carrying out a
statistical procedure with loss function $\mathcal{L}\left( \boldsymbol{C}%
\right) $, then a version of the procedure that minimizes $\mathcal{L}\left(
\eta ^{2}\boldsymbol{I}_{n}\right) $ is \textit{minimax} as $\boldsymbol{C}$
varies over $\mathcal{C}_{M}$. By (ii) this remains true for subsets of $%
\mathcal{C}_{M}$ that contain $\eta ^{2}\boldsymbol{I}_{n}$. To apply this
result we need only ensure that $\eta ^{2}$ is large enough that $\mathcal{C}%
_{M}$ contains the departures, from independence or homoscedasticity, that
are of interest. If so, the designs and estimates of this article -- i.e.
those derived under an assumption of i.i.d.\ errors -- enjoy the additional
optimality property of minimizing the maximum value of $\mathcal{L}\left( 
\boldsymbol{C}\right) $, as the covariance matrix $\boldsymbol{C}$ of $%
\boldsymbol{u}$ varies over $\mathcal{C}_{M}$.

Lemma \ref{lemma: dependence} is unsatisfactory in this M-estimation context
since it is based on the covariance structure of $\left\{ \psi _{\sigma
}\left( \varepsilon _{i}\right) \right\} $, rather than that of $\left\{
\varepsilon _{i}\right\} $. In general there is no tractable relationship
between the two. There are however important exceptions. Suppose that under
the dependence structure the $\left\{ \varepsilon _{i}\right\} $ are \textit{%
exchangeable}, hence \textit{equicorrelated}. Then the $\left\{ \psi
_{\sigma }\left( \varepsilon _{i}\right) \right\} $ also have this
structure: for some $\rho $ (typically \textit{not} the correlation among
the $\left\{ \varepsilon _{i}\right\} $), and some $\alpha ^{2}$, 
\begin{equation*}
\text{\textsc{cov}}\left[ \boldsymbol{u}\right] =\boldsymbol{C}=\alpha
^{2}\left( \left( 1-\rho \right) \boldsymbol{I}_{n}+\rho \boldsymbol{1}_{n}%
\boldsymbol{1}_{n}^{\T}\right) .
\end{equation*}%
We impose bounds $\left\vert \rho \right\vert \leq $ $\rho _{\max }<1$ and $%
\alpha ^{2}\leq \alpha _{\max }^{2}$. We shall work with the norm $%
\left\Vert \boldsymbol{C}\right\Vert _{1}=\alpha ^{2}\left( 1+\left(
n-1\right) \left\vert \rho \right\vert \right) $, this is also $\left\Vert 
\boldsymbol{C}\right\Vert _{E}$ if $\rho \geq 0$. Set $\eta ^{2}=\alpha
_{\max }^{2}\left( 1+\left( n-1\right) \rho _{\max }\right) .$ Then $%
\left\Vert \boldsymbol{C}\right\Vert _{1}\leq \eta ^{2}$ and Lemma \ref%
{lemma: dependence} applies: an assumption of i.i.d.\ errors with \textsc{var}%
$\left[ U_i) \right] =$ $\eta ^{2}$, which for convenience we write as $\eta
^{2}=\eta _{0}^{2}\sigma _{M}^{2}$, is minimax within a class containing all
such equicorrelated error structures.

By a similar treatment one can show robustness against dependence structures
under which $\varepsilon_{i}$ and $\varepsilon_{j}$ are independent for $%
|i-j| > q$, as is the case for $MA(q)$ errors. The case of independent but
heteroscedastic errors is clearly covered as well, by $\alpha _{\max
}^{2}\geq \max_{i}$\textsc{var}$\left[ \psi _{\sigma }\left( \varepsilon
_{i}\right) \right] $.

\section{ Construction of minimax \footnotesize{AIMSE} designs\label{section: construction}}

With 
\begin{equation}
\nu \overset{def}{=}\kappa ^{2}/\left( \eta _{0}^{2}\sigma _{M}^{2}+\kappa
^{2}\right) ,  \label{nu}
\end{equation}%
the results of \S \ref{section: minimax theory} easily accommodate those of 
\S \ref{section: dependence}: the maximum, over both $\tau \in \Upsilon $
and $\boldsymbol{C}\in \mathcal{C}_{1}$, of (\ref{IMSEn}) is given by $%
\left( \eta _{0}^{2}\sigma _{M}^{2}+\kappa ^{2}\right) /n$ (which does not
depend on the design) times%
\begin{equation}
I_{\nu }\left( \xi _{n}\right) =\left( 1-\nu \right) tr\boldsymbol{R}%
^{-1}\left( \xi _{n}\right) +\nu ch_{\max }\boldsymbol{U}\left( \xi
_{n}\right) .  \label{Inu}
\end{equation}%
If robustness against dependence is not of interest then set $\eta _{0}^{2}$ 
$=1$.

A \textit{minimax} \textsc{aimse} design is a minimizer of $I_{\nu }\left( \xi _{n}\right) $%
. For fixed $\nu $, $I_{\nu }\left( \xi _{n}\right) $ is precisely the value
minimized, in \cite{w}, to obtain minimax designs for \textsc{ls} estimates,
thus justifying our statement in \S 1 that such designs remain minimax
optimal, asymptotically, for M-estimation.

The minimization of $I_{\nu }\left( \xi _{n}\right) $ is carried out
sequentially, as described in Theorem 5 of \cite{w}. Briefly, given a
current $k$-point design $\xi _{k}$, the loss resulting from the addition of
a design point at $\boldsymbol{x}_{i}$ is expanded as 
\begin{equation}
I_{\nu }\left( \xi _{k+1}^{(i)}\right) =I_{\nu }\left( \xi _{k}\right)
-t_{k,i}/k+O\left( k^{-2}\right) ,  \label{expansion}
\end{equation}%
and then $\boldsymbol{x}_{(i)}$, with $\left( i\right) =\arg \max_{i}t_{k,i}$%
, is added to the design. This is carried out to convergence, yielding a
design $\xi $ on $\mathcal{\chi }$ with intended allocations $n_{i}=n\xi
\left( \boldsymbol{x}_{i}\right) $. Typically most $\xi \left( \boldsymbol{x}%
_{i}\right) $ are zero, but otherwise the $n_{i}$ are not integers. To
obtain implementable designs $\xi_{n}$ we first round up the $n_{i}$ to $%
\left\lceil n\xi \left( \boldsymbol{x}_{i}\right) \right\rceil $, whose sum
exceeds $n$. The excess is decreased stepwise, by removing points whose
value of $t_{k,i}$ in (\ref{expansion}) is a minimum. This method typically
results in only a very small increase in the minimized value of $I_{\nu
}\left( \xi _{n}\right) $.

\section{Integrated M-loss\label{section: IMloss}}

Our method of measuring loss -- \textsc{imse} -- is itself heavily
influenced by outlying responses, and hence is non-robust. As an
alternative, in this section we consider \textit{Integrated M-loss }(\textsc{%
im}), in which \textsc{imse} is replaced by 
\begin{equation*}
\text{\textsc{im}}\left( \xi ,\tau \right) =\sum_{\boldsymbol{x\in \chi }}E%
\left[ \rho \left( E\left[ Y\left( \boldsymbol{x}\right) \right] -\hat{Y}%
\left( \boldsymbol{x}\right) \right) \right] =\sum_{\boldsymbol{x\in \chi }}E%
\left[ \rho \left( \tau \left( \boldsymbol{x}\right) -\boldsymbol{f}^{\T}\left( \boldsymbol{x}\right) \left( 
\boldsymbol{\hat{\theta}}_{n}\left( \xi |\tau \right) -\boldsymbol{\theta }%
_{0}\left( \tau \right) \right) \right) \right] ,
\end{equation*}%
where $\rho \left( x\right) =\int_{0}^{x}\psi \left( t\right) dt$. This is
evaluated at the asymptotic distribution (\ref{AN}), yielding `Asymptotic Integrated M-loss' (\textsc{aim}). If $\psi $ is Huber's
score function then \textsc{im} is termed `Huber loss' (\cite{htf09}, p.\ 349). If $\psi \left( x\right) =x$ then \textsc{im} is identical to 
\textsc{imse} applied to \textsc{ls} estimates. We continue to impose constraints (\ref{orthogonality}) and (%
\ref{bound}) on $\tau \left( \boldsymbol{x}\right) $. We define  
$
\boldsymbol{H}\left( \boldsymbol{\xi }_{n}\right) =\left( \boldsymbol{I}_{N}-%
\boldsymbol{Q\boldsymbol{R}}^{-1}\left( \boldsymbol{\xi }_{n}\right) 
\boldsymbol{Q}^{\T}\boldsymbol{D}\left( \boldsymbol{\xi }%
_{n}\right) \right) \boldsymbol{Q}_{\perp }:N\times N-p
$, 
where $\boldsymbol{\boldsymbol{R}}\left( \boldsymbol{\xi }\right) $ is as in
Theorem \ref{thm: max loss} and, as in the proof of that theorem, $\boldsymbol{Q}_{\boldsymbol{\perp }}:N\times N-p$ is the orthogonal
complement of $\boldsymbol{Q}$.

\begin{lemma}
\label{lemma: max IMloss}
Asymptotic Integrated M-loss is given by 
$\textsc{aim} =\max_{\left\Vert \boldsymbol{\beta}\right\Vert \leq 1} G_{0}\left( \boldsymbol{\beta }\right)$, 
for
\begin{equation}
G_{0}\left( \boldsymbol{\beta }\right) =\sum_{i=1}^{N}\int_{-\infty
}^{\infty }\rho \left( \frac{\kappa }{\sqrt{n}}\boldsymbol{h}_{i}^{\T}\left( \boldsymbol{\xi }_{n}\right) \boldsymbol{\beta 
}-\frac{\sigma _{M}}{\sqrt{n}}s_i\left( \boldsymbol{\xi }_{n}\right)z\right) \phi \left( z\right) dz,
\label{im}
\end{equation}
where $\boldsymbol{h}_{i}^{\T}\left(\boldsymbol{\xi }_{n}\right) $ is the i$^{th}$ row of 
$\boldsymbol{H}\left( \boldsymbol{\xi }_{n}\right)$,
$s_{i}\left( \boldsymbol{\xi }_{n}\right) =\sqrt{\left( \boldsymbol{QR}%
^{-1}\left( \boldsymbol{\xi }_{n}\right) \boldsymbol{Q}^{\T}\right) _{ii}}$, 
and $\phi \left( z\right) $ is the standard Normal density.
\end{lemma}

The computation of minimax \textsc{aim} designs is
expensive, since the maximization of $G_{0}\left( \boldsymbol{\beta }\right) 
$, done analytically for \textsc{imse}, must now be done numerically. For this we employ \textsc{matlab}'s constrained  minimizer \texttt{fmincon} to
minimize $-G_{0}\left( \boldsymbol{\beta }\right) $ over $\left\Vert 
\boldsymbol{\beta }\right\Vert \leq 1$. This is sped up by supplying the
gradient and Hessian, and by vectorizing the input functions. The integrals
are approximated via Gauss-Hermite quadrature. In Lemma \ref{lemma: derivs}
below these integrals are presented in the form 
$I=\int_{-\infty }^{\infty }{g}\left(\boldsymbol{\eta}\left(z%
\right)\right)\phi(z) dz$ for a vectorized function $g$. This is an $N\times 1$ vector with $i^{th}$ element $\int_{-\infty }^{\infty
}g\left( \eta_i\left(z\right)\right) \phi \left( z\right) dz$. This vector of
integrals can be approximated very efficiently as 
$\hat{I}=g\left(\boldsymbol{\eta}\left( \boldsymbol{z}^{\T} \right) \right)\boldsymbol{w}$, 
where the matrix $g\left(\boldsymbol{\eta}\left( \boldsymbol{z}^{\T} \right) \right)$ 
has $\left( i,j\right) ^{th}$ element $g\left(\eta_{i}\left(z_{j}\right)\right) $, for
Gauss-Hermite quadrature nodes $z_{j}$, with corresponding weights $w_{j}$ in $\boldsymbol{w}$.

\begin{lemma}
\label{lemma: derivs} Assume that $%
\rho $, $\psi $ and $\psi ^{\prime }$ are vectorized, so that if $\boldsymbol{M}$ is a
matrix then ${\rho }\left( \boldsymbol{M}\right) $ is a matrix of
the same size with $\left[ {\rho }\left( \boldsymbol{M}\right) %
\right] _{ij}=\rho \left( \boldsymbol{M}_{ij}\right) $, and similarly with $%
{\psi }$ and ${\psi }^{\prime }$. Let 
$\boldsymbol{s}\left( \boldsymbol{\xi }_{n}\right)$ have elements 
$s_i\left( \boldsymbol{\xi }_{n}\right) $
and define
\begin{equation*}
\boldsymbol{\eta}\left( \boldsymbol{\xi }_{n},z%
\right) =  \frac{\kappa }{%
\sqrt{n}}\boldsymbol{H}\left( \boldsymbol{\xi }_{n}\right) \boldsymbol{\beta 
}-\frac{\sigma _{M}}{\sqrt{n}}\boldsymbol{s}\left( \boldsymbol{\xi }%
_{n}\right) z: N \times 1.  
\end{equation*}
Then 
\begin{equation}
G_{0}\left( \boldsymbol{\beta }\right) =\boldsymbol{1}_{N}^{\T}\int_{-\infty }^{\infty }{\rho }
\left(\boldsymbol{\eta}\left( \boldsymbol{\xi }_{n},z\right) \right)  \phi \left( z\right) dz,  \label{G}
\end{equation}%
with $N-p\times 1$ gradient%
\begin{equation}
\boldsymbol{G}_{1}\left( \boldsymbol{\beta }\right) =\frac{\kappa }{\sqrt{n}}%
\boldsymbol{H}^{\T}\left( \boldsymbol{\xi }%
_{n}\right) \int_{-\infty }^{\infty }\psi\left(\boldsymbol{\eta}\left( \boldsymbol{\xi }_{n},z\right) \right)  \phi \left( z\right) dz,  \label{G1}
\end{equation}
and $N-p\times N-p$ Hessian 
\begin{equation}
\boldsymbol{G}_{2}\left( \boldsymbol{\beta }\right) =\frac{\kappa ^{2}}{n}%
\boldsymbol{H}^{\T}\left( \boldsymbol{\xi }%
_{n}\right)
\oplus _{i=1}^{N}  \Big[  \int_{-\infty }^{\infty }\psi^{\prime}\left(\boldsymbol{\eta}\left( \boldsymbol{\xi }_{n},z\right) \right)\phi \left( z\right) dz \Big]_{i}
 \boldsymbol{H}\left( \xi_{n}\right) .  \label{G2}
\end{equation}
\end{lemma}

The subsequent minimization of $\max_{\left\Vert \boldsymbol{\beta}\right\Vert \leq 1} G_{0}\left( \boldsymbol{\beta }\right)$ over $\boldsymbol{\xi }_{n}$ is also carried out by \texttt{fmincon}; this
is quite time-consuming since each trial design $\boldsymbol{\xi }_{n}$ must be
accompanied by a maximizing $\boldsymbol{\beta }$. This latter point renders the sequential approach to the minimization, as employed in  \S \ref{section: construction}, untenable due to the huge number of maximizations required. The final design returned by \texttt{fmincon} is made implementable by applying the `Efficient Design Apportionment' method of \cite{pr92}.

\section{Examples, implementations and discussion \label{section: examples}}

\subsection{Minimax \textsc{aimse} designs} \label{section: imse examples}

To obtain the minimax \textsc{aimse} design only $\nu $, at (\ref{Inu}), need be specified
by the user. The simplest and most natural way to do this is to view $\nu $
as expressing the emphasis on the reduction of losses due to bias, rather
than to variation. Its choice is then up to the user; quite typically $\nu
=.5$ is chosen. In this method there is \textit{no} difference in the minimax designs
for different M-estimates.

Another method is suggested by the definition of $\nu $ at (\ref{nu}).
Although the parameters involved in this definition would not be known to
the designer, it is of interest to see how their values could affect the
resulting designs.

We begin by seeing how much $\nu $ can change from its value under \textsc{ls%
}. The relationship between $\nu $ using Least Squares and $\nu $ using the
M-estimate is, with $\gamma \overset{def}{=}\eta _{0}\sigma /\kappa $, that 
\begin{equation*}
\nu _{\text{\textsc{ls}}}=\left( \gamma ^{2}+1\right) ^{-1}\text{ and }\nu _{%
\text{\textsc{m}}}=\left( \frac{\gamma ^{2}\sigma _{M}^{2}}{\sigma ^{2}}%
+1\right) ^{-1}.
\end{equation*}%
We assess these assuming that $\varepsilon /\sigma \sim N(0,1)$ and that $%
\psi \left( x\right) =xI\left( \left\vert x\right\vert \leq c\right)
+sign(x)cI\left( \left\vert x\right\vert >c\right) $ (\cite{h64}).

\begin{lemma}
\label{lemma: maxdiff} With notation as above, 
\begin{equation*}
0\leq \nu _{\text{\textsc{ls}}}-\nu _{\text{\textsc{m}}}\leq \frac{\sqrt{\pi
/2}-1}{\sqrt{\pi /2}+1}\simeq .1124.
\end{equation*}%
The lower bound is attained only when the M-estimate is the \textsc{lse},
and the upper bound is attained when the M-estimate is the $L_{1}$ estimate
and $\gamma ^{2}=1/\sqrt{\pi /2}\simeq .7979$. At the maximum 
\begin{equation}
\nu _{\text{\textsc{ls}}}=\frac{\sqrt{\pi /2}}{\sqrt{\pi /2}+1}\simeq
.5562,\quad \nu _{\text{\textsc{m}}}=\frac{1}{\sqrt{\pi /2}+1}=1-\nu _{\text{%
\textsc{ls}}}\simeq .4438.  \label{nus}
\end{equation}
\end{lemma}

\begin{figure}[tbp]
\centering
\includegraphics[scale=1.2]{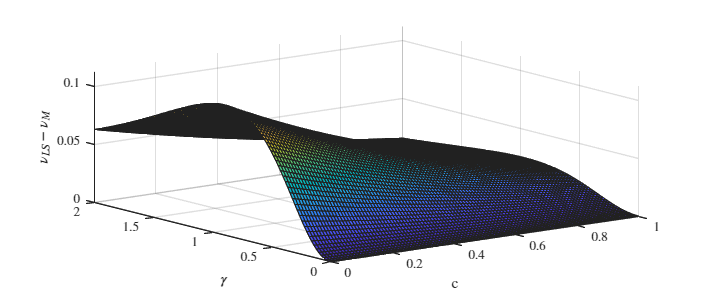}
\caption{Differences $\protect\nu _{\text{\textsc{ls}}}-\protect\nu _{\text{%
\textsc{m}}}$ in terms of $\protect\gamma $ and $c$. }
\label{fig:diff}
\end{figure}

\begin{figure}[tbp]
\centering
\includegraphics[scale=1.2]{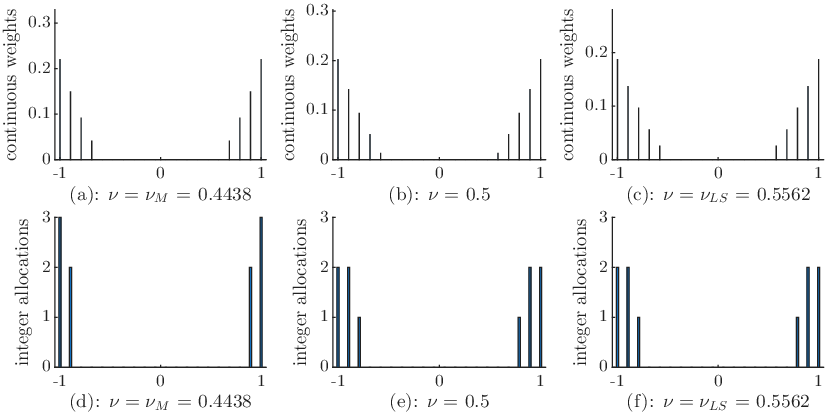}
\caption{Designs for linear regression; n=10, N=20. (a)-(c): Continuous
weights; (d)-(f): Integer allocations.}
\label{fig:lindes}
\end{figure}

\begin{figure}[tbp]
\centering
\includegraphics[scale=1.2]{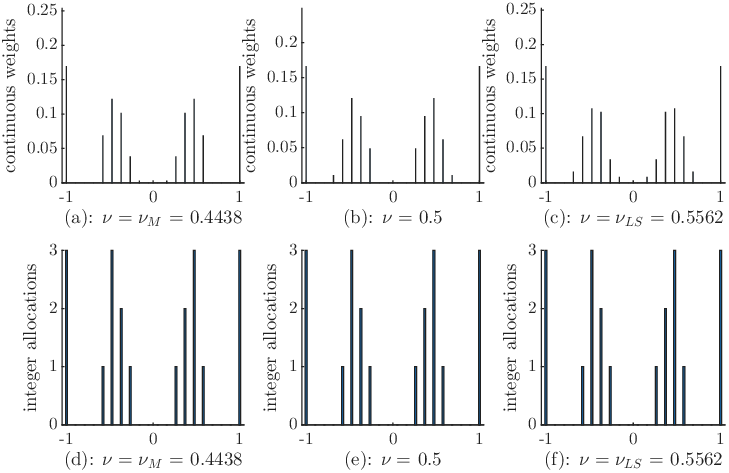}
\caption{Designs for cubic regression; n=20, N=20. (a)-(c): Continuous
weights; (d)-(f): Integer allocations.}
\label{fig:cubdes}
\end{figure}

See Figure \ref{fig:diff}. We have constructed designs, using the
`worst case' values $\nu _{\text{\textsc{ls}}}$ and $\nu _{\text{\textsc{m}}}$ of (\ref{nus}), and also their
midpoint $\nu =.5$. See Figure \ref{fig:lindes} for approximate linear
regression and Figure \ref{fig:cubdes} for approximate cubic regression,
each on design spaces consisting of 20 equally spaced points spanning $[-1,1]
$. Given that $\nu _{\text{\textsc{ls}}}$, $\nu _{\text{\textsc{m}}}$ and
the M-estimate were chosen so as to make the designs as different as
possible, these designs are remarkably similar -- even identical after being
made implementable. The implementations all can be described as taking the
replicates that would otherwise be assigned by the classically ($\nu =0$)
I-optimal designs and spreading them out into clusters at nearby design
points. The I-optimal design for linear regression places mass of $.5$ at
each of $\pm 1$. That for cubic regression was derived by \cite{s77} and
places masses of $.1545$ and $.3455$ at $\pm 1$ and $\pm .4472$.

We conclude that an experimenter measuring loss by \textsc{imse} should feel quite safe in using the same
design for an experiment regardless of which M-estimate is to be employed,
and in choosing $\nu $ to represent his desired emphasis on bias reduction,
as posited at the beginning of this section. As well, the assumption of
i.i.d. errors is safe, even optimal in a minimax sense, amongst
equicorrelated or heteroscedastic error structures (and others).

\subsection{Minimax \textsc{aim} designs} \label{section: im examples}

\begin{figure}[tbp]
\centering
\includegraphics[scale=1.2]{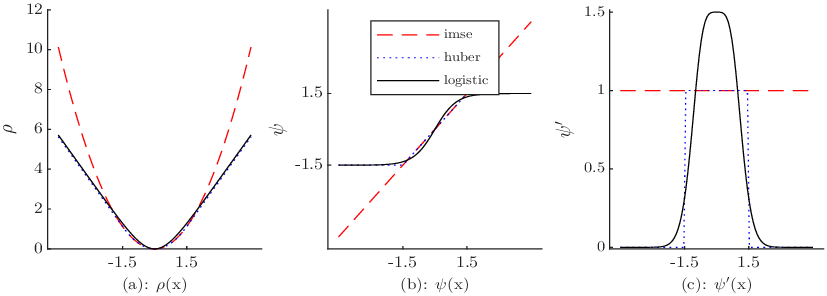}
\caption{Loss functions $\rho(x)$ with their derivatives for \textsc{im} designs. \textsc{imse} loss: $\psi(x) = x$; Huber loss: $\psi \left( x\right) =xI\left( \left\vert x\right\vert \leq 1.5\right) +sign(x)1.5I\left( \left\vert x\right\vert >1.5\right) $; logistic loss: $\psi(x) = 1.5\tanh(x)$.}  
\label{fig:scores}
\end{figure}

\begin{figure}[tbp]
\centering
\includegraphics[scale=1.2]{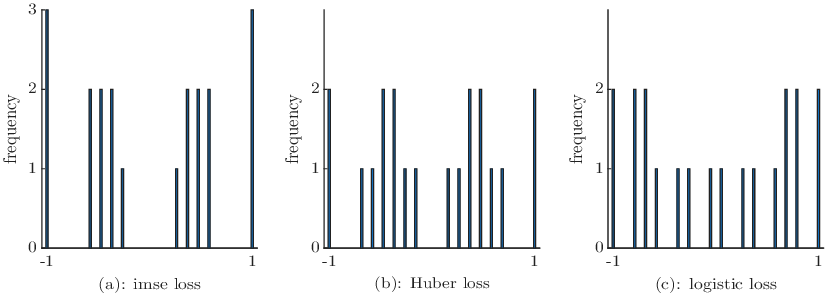}
\caption{Minimax \textsc{aim} designs for cubic regression; n=20, N=20.}  
\label{fig:im1}
\end{figure}

\begin{figure}[tbp]
\centering
\includegraphics[scale=1.2]{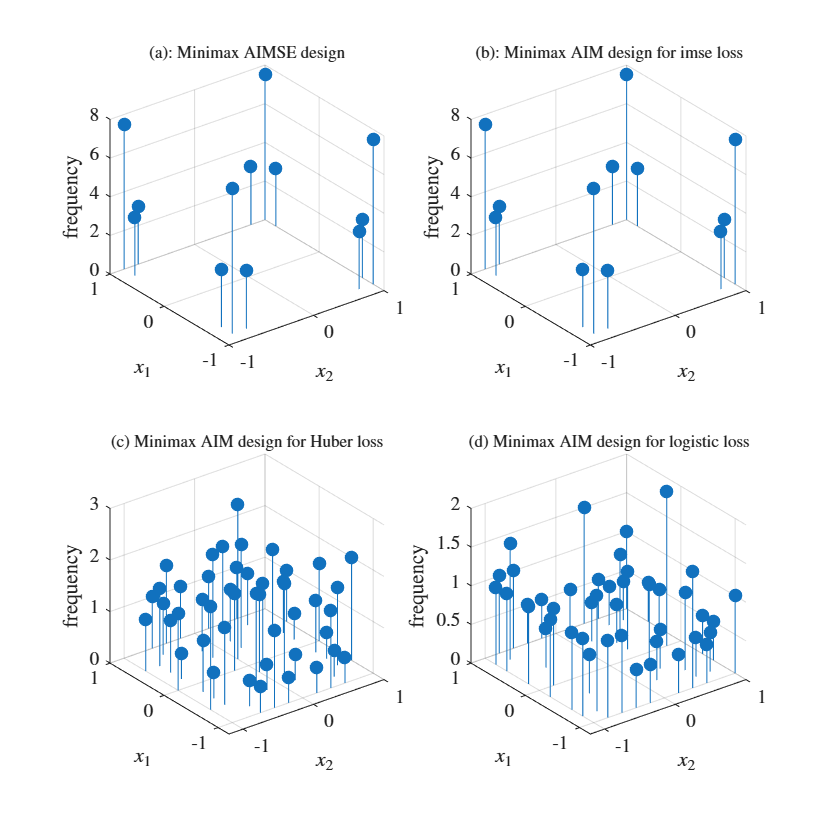}
\caption{Minimax \textsc{aim} designs for bivariate linear regression; n=36, N=121.}  
\label{fig:bivar}
\end{figure}

We have applied the procedure described in \S \ref{section: IMloss} using the three loss functions illustrated in Figure \ref {fig:scores}. Figure  \ref {fig:im1} gives the results for approximate cubic regression, using the same values of $n$ and $N$ as in Figure \ref {fig:cubdes}, and $\kappa = \sigma_M = 5$ (hence $\nu = .5$). The function being minimized has local minima, and we had to run the algorithm repeatedly, using a different starting design to initialize the iterative process each time, before arriving at what appears to be at least very close to an absolute minimum. The designs are qualitatively very similar to those of \ref {fig:cubdes}(e), the differences in the clusters arising from the different implementation method being used, as well as the robust loss functions used in Figures \ref {fig:im1}(b) and \ref {fig:im1}(c), which typically result in more diffuse designs. This latter feature is quite pronounced in Figure \ref {fig:bivar}, in which the fitted model is linear in two variables $x_1, x_2\in \left[-1,1\right]$. The designs in Figures \ref {fig:bivar}(a) and \ref {fig:bivar}(b) are based on squared error loss and, as expected, have clusters of design points in the corners of the design space. The robustness of the designs in Figures \ref {fig:bivar}(c) and \ref{fig:bivar}(d) is reflected in the increased spread of these clusters.

\appendix

\section*{Appendix: Proofs}

\setcounter{equation}{0} \renewcommand{\theequation}{A.\arabic{equation}} %
\setcounter{lemma}{0}
\renewcommand{\thelemma}{A.\arabic{lemma}} 
\setcounter{subsection}{0} \renewcommand{\thesubsection}{A.%
\arabic{subsection}} 

\begin{remark}
It is convenient to identify a design with its design measure -- a discrete
probability measure $\xi _{n}\left( d\boldsymbol{x}\right) $ on $\mathcal{X}$%
. Thus, as an example, if the design calls for $n_{i}$ runs to be made at $%
\boldsymbol{x}_{i}\in \mathcal{\chi }$, we write averages of functions $\phi 
$ of the data as 
\begin{equation*}
\frac{1}{n}\sum_{\left\{ i|n_{i}>0\right\} }\sum_{j=1}^{n_{i}}\phi \left(
Y_{j}\left( \boldsymbol{x}_{i}\right) \right) =\sum_{\left\{
i|n_{i}>0\right\} }\frac{n_{i}}{n}\bar{\phi}\left( \boldsymbol{x}_{i}\right)
=\int_{\mathcal{X}}\bar{\phi}\left( \boldsymbol{x}\right) \xi _{n}\left( d%
\boldsymbol{x}\right) ,
\end{equation*}%
where $\bar{\phi}\left( \boldsymbol{x}_{i}\right) =\left( 1/n_{i}\right)
\sum_{j=1}^{n_{i}}\phi \left( Y_{j}\left( \boldsymbol{x}_{i}\right) \right) $%
. We abbreviate $\sum_{\left\{ i|n_{i}>0\right\} }\sum_{j=1}^{n_{i}}$ by $%
\sum_{i,j}$.
\end{remark}

Applying the theorem of \cite{vw98} mentioned in Remark \ref{remark:
assumptions} requires that we verify further properties of the limit, given
in Lemma \ref{lemma: limit} below.

\begin{lemma}
\label{lemma: limit} As well as C1) - C3) assume that the parameter space $%
\Theta $ is a compact subset of $\mathbb{R}^{p}$. Then with $\Psi _{n}\left( 
\boldsymbol{\theta }\right) \overset{def}{=}n^{-1}\sum_{i,j}\psi _{\sigma
}\left( Y_{j}\left( \boldsymbol{x}_{i}\right) -\boldsymbol{f}^{\T}\left( \boldsymbol{x}_{i}\right) \left( \boldsymbol{%
\theta }+\boldsymbol{\tilde{\theta}}_{n}\right) \right) \boldsymbol{f}\left( 
\boldsymbol{x}_{i}\right) ,$ the function 
\begin{equation*}
\Psi \left( \boldsymbol{\theta }\right) =\int_{\mathcal{\chi }}E\left[ \psi
_{\sigma }\left( Y\left( \boldsymbol{x}\right) -\boldsymbol{f}^{\T}\left( \boldsymbol{x}\right) \left( \boldsymbol{%
\theta }+\boldsymbol{\theta }_{0}\right) \right) \right] \boldsymbol{f}%
\left( \boldsymbol{x}\right) \xi _{\ast }\left( d\boldsymbol{x}\right)
\end{equation*}%
satisfies 
\begin{eqnarray*}
(i) &&\sup_{\boldsymbol{\theta }\in \boldsymbol{\Theta }}\left\Vert \Psi
_{n}\left( \boldsymbol{\theta }\right) -\Psi \left( \boldsymbol{\theta }%
\right) \right\Vert \overset{pr}{\rightarrow }0, \\
(ii) &\text{ }&\text{for every }\delta >0\text{, }\inf_{\left\Vert 
\boldsymbol{\theta }\right\Vert \geq \delta }\left\Vert \Psi \left( 
\boldsymbol{\theta }\right) \right\Vert >0=\left\Vert\Psi ( \boldsymbol{0}%
_p\right)\Vert .
\end{eqnarray*}
\end{lemma}

\noindent \textit{Proof}: With $\boldsymbol{\tilde{\Psi}}_{n}(\boldsymbol{\theta })\overset{def}{=}%
\int_{\chi }E\!\left[ \psi _{\sigma }\!\left( Y(\boldsymbol{x})-\boldsymbol{f%
}^{\T}(\boldsymbol{x})(\boldsymbol{\theta }+%
\widetilde{\boldsymbol{\theta }}_{n})\right) \right] \boldsymbol{f}(%
\boldsymbol{x})\,\xi _{n}(d\boldsymbol{x})=E\left[ \Psi _{n}\left( 
\boldsymbol{\theta }\right) \right] $ we shall establish (i) of the Lemma by
showing:

\begin{enumerate}
\item[(a)] $\sup_{\boldsymbol{\theta }\in \boldsymbol{\Theta }}\left\Vert 
\boldsymbol{\tilde{\Psi}}_{n}\left( \boldsymbol{\theta }\right) -\boldsymbol{%
\Psi }\left( \boldsymbol{\theta }\right) \right\Vert \rightarrow 0$,

\item[(b)] $\sup_{\boldsymbol{\theta }\in \boldsymbol{\Theta }}\left\Vert 
\boldsymbol{\Psi }_{n}\left( \boldsymbol{\theta }\right) -\boldsymbol{\tilde{%
\Psi}}_{n}\left( \boldsymbol{\theta }\right) \right\Vert \overset{pr}{%
\rightarrow }0$.
\end{enumerate}

For (a), with $\boldsymbol{h}(\boldsymbol{x},\boldsymbol{\theta },%
\boldsymbol{\eta })\overset{def}{=}E\!\left[ \psi _{\sigma }\!\left( Y(%
\boldsymbol{x})-\boldsymbol{f}^{\T}(%
\boldsymbol{x})(\boldsymbol{\theta }+\boldsymbol{\eta })\right) \right] 
\boldsymbol{f}(\boldsymbol{x})$ we have $\boldsymbol{\Psi }(\boldsymbol{%
\theta })=\int_{\chi }\boldsymbol{h}(\boldsymbol{x},\boldsymbol{\theta },%
\boldsymbol{\theta }_{0})\,\xi _{\ast }(d\boldsymbol{x})$ and $\boldsymbol{%
\tilde{\Psi}}_{n}\left( \boldsymbol{\theta }\right) =\int_{\chi }\boldsymbol{%
h}(\boldsymbol{x},\boldsymbol{\theta },\boldsymbol{\tilde{\theta}}_{n})\,\xi
_{n}(d\boldsymbol{x})$, with 
\begin{eqnarray*}
\boldsymbol{\Psi }(\boldsymbol{\theta })-\boldsymbol{\tilde{\Psi}}_{n}\left( 
\boldsymbol{\theta }\right) &=&\int_{\chi }\left( \boldsymbol{h}(\boldsymbol{%
x},\boldsymbol{\theta },\boldsymbol{\theta }_{0})-\boldsymbol{h}(\boldsymbol{%
x},\boldsymbol{\theta },\boldsymbol{\tilde{\theta}}_{n})\right) \,\xi _{n}(d%
\boldsymbol{x})+\int_{\chi }\boldsymbol{h}(\boldsymbol{x},\boldsymbol{\theta 
},\boldsymbol{\theta }_{0})\,(\xi _{\ast }-\xi _{n})(d\boldsymbol{x}) \\
&=&\boldsymbol{\alpha }_{n}(\boldsymbol{\theta })+\boldsymbol{\beta }_{n}(%
\boldsymbol{\theta })\text{, say.}
\end{eqnarray*}%
Since $\boldsymbol{h}(\boldsymbol{x},\boldsymbol{\theta },\boldsymbol{\eta }%
) $ is differentiable w.r.t. $\boldsymbol{\eta }$, and both $\psi _{\sigma
}^{\prime }$ and $\boldsymbol{f}$ are bounded, there is $C>0$ for which $%
\left\Vert \boldsymbol{h}(\boldsymbol{x},\boldsymbol{\theta },\boldsymbol{%
\theta }_{0})-\boldsymbol{h}(\boldsymbol{x},\boldsymbol{\theta },\boldsymbol{%
\tilde{\theta}}_{n})\right\Vert \leq C\left\Vert \boldsymbol{\theta }_{0}-%
\boldsymbol{\tilde{\theta}}_{n}\right\Vert $, uniformly in $\boldsymbol{x}%
\in \chi $ and $\boldsymbol{\theta }\in \Theta $. Hence 
\begin{equation}
\sup_{\boldsymbol{\theta }\in \Theta }\left\Vert \boldsymbol{\alpha }_{n}(%
\boldsymbol{\theta })\right\Vert \leq C\left\Vert \boldsymbol{\theta }_{0}-%
\boldsymbol{\tilde{\theta}}_{n}\right\Vert \rightarrow 0.  \label{first sup}
\end{equation}%
With $\boldsymbol{g}_{\boldsymbol{\theta }}(\boldsymbol{x})\overset{def}{=}%
\boldsymbol{h}(\boldsymbol{x},\boldsymbol{\theta },\boldsymbol{\theta }_{0})$
we obtain, in a similar fashion, $\sup_{\boldsymbol{x}\in \chi }\Vert 
\boldsymbol{g}_{\boldsymbol{\theta }_{1}}(\boldsymbol{x})-\boldsymbol{g}_{%
\boldsymbol{\theta }_{2}}(\boldsymbol{x})\Vert \leq C\Vert \boldsymbol{%
\theta }_{1}-\boldsymbol{\theta }_{2}\Vert $. Then since $\boldsymbol{\Theta 
}$ is compact, for every $\delta >0$ there exist $\boldsymbol{\theta }%
^{(1)},\dots ,\boldsymbol{\theta }^{(K)}\in \Theta $ such that every $%
\boldsymbol{\theta }\in \boldsymbol{\Theta }$ lies within distance $\delta $
of one of these points. Therefore, 
\begin{equation*}
\sup_{\boldsymbol{\theta }\in \Theta }\left\Vert \boldsymbol{\beta }_{n}(%
\boldsymbol{\theta })\right\Vert \leq 2C\delta +\max_{1\leq k\leq
K}\left\Vert \int_{\chi }\boldsymbol{g}_{\boldsymbol{\theta }^{(k)}}(%
\boldsymbol{x})\,\xi _{\ast }(d\boldsymbol{x})-\int_{\chi }\boldsymbol{g}_{%
\boldsymbol{\theta }^{(k)}}(\boldsymbol{x})\,\xi _{n}(d\boldsymbol{x}%
)\right\Vert .
\end{equation*}%
That $\xi _{n}\overset{d}{\rightarrow }\xi _{\ast }$ implies that the
preceding norms, hence their maximum, converge to $0$. Letting $\delta
\rightarrow 0$, we conclude that $\sup_{\boldsymbol{\theta }\in \Theta
}\left\Vert \boldsymbol{\beta }_{n}(\boldsymbol{\theta })\right\Vert
\rightarrow 0$, which together with (\ref{first sup}) yields (a).

For (b), we write $\boldsymbol{\Delta }_{n}(\boldsymbol{\theta })\overset{def%
}{=}\boldsymbol{\Psi }_{n}(\boldsymbol{\theta })-\boldsymbol{\tilde{\Psi}}%
_{n}\left( \boldsymbol{\theta }\right) \ =n^{-1}\sum_{i,j}\left\{ 
\boldsymbol{\zeta }_{ij,n}(\boldsymbol{\theta })-E[\boldsymbol{\zeta }%
_{ij,n}(\boldsymbol{\theta })]\right\} ,$ where 
\begin{equation*}
\boldsymbol{\zeta }_{ij,n}(\boldsymbol{\theta })=\psi _{\sigma }\!\left(
Y_{j}(\boldsymbol{x}_{i})-\boldsymbol{f}^{\T}
\left( \boldsymbol{x}_{i}\right) (\boldsymbol{\theta }+\widetilde{%
\boldsymbol{\theta }}_{n})\right) \boldsymbol{f}(\boldsymbol{x}_{i}).
\end{equation*}%
Arguing as above, there exists $C>0$ such that for all $\boldsymbol{\theta }%
_{1},\boldsymbol{\theta }_{2}\in \boldsymbol{\Theta }$ and each pair $(i,j)$%
, 
\begin{equation*}
\left\Vert \boldsymbol{\zeta }_{ij,n}(\boldsymbol{\theta }_{1})-\boldsymbol{%
\zeta }_{ij,n}(\boldsymbol{\theta }_{2})\right\Vert \leq C\left\Vert 
\boldsymbol{\theta }_{1}-\boldsymbol{\theta }_{2}\right\Vert .
\end{equation*}%
Taking expectations gives the same bound for the centred version, so that 
\begin{equation*}
\left\Vert \boldsymbol{\Delta }_{n}(\boldsymbol{\theta }_{1})-\boldsymbol{%
\Delta }_{n}(\boldsymbol{\theta }_{2})\right\Vert \leq 2C\left\Vert 
\boldsymbol{\theta }_{1}-\boldsymbol{\theta }_{2}\right\Vert .
\end{equation*}%
Thus $\Delta _{n}$ is uniformly Lipschitz on $\boldsymbol{\Theta }$, with a
constant independent of $n$.

Fix $\varepsilon >0$. Choose $\delta >0$ such that $2C\delta <\varepsilon /2$%
, and let $\boldsymbol{\theta }^{(1)},\dots ,\boldsymbol{\theta }^{(K)}$ be
a finite $\delta $-net of $\Theta $, as above. Then 
\begin{equation*}
\sup_{\boldsymbol{\theta }\in \Theta }\left\Vert \boldsymbol{\Delta }_{n}(%
\boldsymbol{\theta })\right\Vert \leq \max_{1\leq k\leq K}\left\Vert 
\boldsymbol{\Delta }_{n}(\boldsymbol{\theta }^{(k)})\right\Vert +\varepsilon
/2,
\end{equation*}%
and it remains only to show that 
\begin{equation}
\max_{1\leq k\leq K}\left\Vert \boldsymbol{\Delta }_{n}(\boldsymbol{\theta }%
^{(k)})\right\Vert \overset{pr}{\rightarrow }0.  \label{remains}
\end{equation}%
Write $\Delta _{n}(\boldsymbol{\theta }^{(k)})=n^{-1}\sum_{i,j}\left( 
\boldsymbol{\zeta }_{ij,n}(\boldsymbol{\theta }^{(k)})-E\left[ \boldsymbol{%
\zeta }_{ij,n}(\boldsymbol{\theta }^{(k)})\right] \right) =n^{-1}\sum_{i,j}%
\boldsymbol{Z}_{ij,n}^{(k)}$, say. We can write $\boldsymbol{\zeta }_{ij,n}(%
\boldsymbol{\theta })$ as $\boldsymbol{\zeta }_{ij,n}(\boldsymbol{\theta }%
)=\psi _{\sigma }\!\left( \varepsilon _{ij}+e_{i,n}\left( \boldsymbol{\theta 
}\right) \right) \boldsymbol{f}(\boldsymbol{x}_{i}) $, for certain functions 
$e_{i,n}\left( \boldsymbol{\theta }\right) $, continuous on $\boldsymbol{%
\Theta }$ hence uniformly bounded there. Then $E[\boldsymbol{Z}%
_{ij,n}^{(k)}]=\boldsymbol{0}$ and{\ there is }$C>0$ for which {$%
\sup_{i,j,n,k}E\Vert \boldsymbol{Z}_{ij,n}^{(k)}\Vert ^{2}\leq C$}. Since
the $\boldsymbol{Z}_{ij,n}^{(k)}$ are independent for each $(k,n)$, \textsc{%
var}$\left[ \Delta _{n}(\boldsymbol{\theta }^{(k)})\right] \leq
n^{-1}C\rightarrow 0$, implying by Chebyshev's inequality that $\boldsymbol{%
\Delta }_{n}(\boldsymbol{\theta }^{(k)})\overset{pr}{\rightarrow }0$ for
each $k=1,\dots ,K$. Now (\ref{remains}) follows, thus completing the proof
of (b).

For part (ii) of the Lemma, recall that the design space is finite and that
each $\xi _{n}\left( \boldsymbol{x}_{i}\right) $ tends to $\xi _{\ast
}\left( \boldsymbol{x}_{i}\right) $, so that 
\begin{equation}
\int_{\mathcal{\chi }}E\left[ \psi _{\sigma }\left( Y\left( \boldsymbol{x}%
\right) -\boldsymbol{f}^{\T}\left( 
\boldsymbol{x}\right) \boldsymbol{\theta }_{0}\right) \right] \boldsymbol{f}%
\left( \boldsymbol{x}\right) \xi _{\ast }\left( d\boldsymbol{x}\right)
=\lim_{n\rightarrow \infty }\int_{\mathcal{\chi }}E\left[ \psi _{\sigma
}\left( \varepsilon +\tau \left( \boldsymbol{x}\right) \right) \right] 
\boldsymbol{f}\left( \boldsymbol{x}\right) \xi _{n}\left( d\boldsymbol{x}%
\right) =\boldsymbol{0}_p,  \label{root}
\end{equation}%
since $\tau \left( \boldsymbol{x}\right) $ is $O\left( n^{-1/2}\right) $. It
follows that $\Psi \left( \boldsymbol{0}_p\right) =\boldsymbol{0}_p $.
By (\ref{root}), $\boldsymbol{\theta }_{0}$ is a stationary point of the
function $k\left( \boldsymbol{\theta }\right) \overset{def}{=}\int_{\mathcal{%
\chi }}E\left[ \rho _{\sigma }\left( Y\left( \boldsymbol{x}\right) -%
\boldsymbol{f}^{\T}\left( \boldsymbol{x}%
\right) \boldsymbol{\theta }\right) \right] \xi _{\ast }\left( d\boldsymbol{x%
}\right) $, where $\rho _{\sigma }(x)=\int_{0}^{x}\psi _{\sigma }\left( t\right)
dt$. The Hessian is $\ddot{k}\left( \boldsymbol{\theta }\right) =\int_{%
\mathcal{\chi }}E\left[ \psi _{\sigma }^{\prime }\left( \varepsilon +M_{%
\boldsymbol{\theta }}\left( \boldsymbol{x}\right) \right) \right] 
\boldsymbol{f}\left( \boldsymbol{x}\right) \boldsymbol{f}^{\T}\left( \boldsymbol{x}\right) \xi _{\ast }\left( d\boldsymbol{x}\right)$ , 
for $M_{\boldsymbol{\theta }}\left( \boldsymbol{x}%
\right) =E\left[ Y\left( \boldsymbol{x}\right)] -\boldsymbol{f}^{\T}\left( \boldsymbol{x}\right) \boldsymbol{\theta }%
\right) $. Since $\psi_{\sigma }$ is weakly increasing, $\psi _{\sigma
}^{\prime }\left( \varepsilon +M_{\boldsymbol{\theta }}\left( \boldsymbol{x}%
\right) \right) \geq 0$ and so $E\left[ \psi _{\sigma }^{\prime }\left(
\varepsilon +M_{\boldsymbol{\theta }}\left( \boldsymbol{x}\right) \right) %
\right] $ must be strictly positive, else $\psi _{\sigma }^{\prime }\left(
\varepsilon +M_{\boldsymbol{\theta }}\left( \boldsymbol{x}\right) \right)
\equiv 0$ (a.s.).
 Put 
\begin{equation*}
c\left( \boldsymbol{\theta }\right) =\min_{i=1,...,N}\left\{ E\left[ \psi
_{\sigma }^{\prime }\left( \varepsilon +M_{\boldsymbol{\theta }}\left( 
\boldsymbol{x}_{i}\right) \right) \right] \left\vert {}\right. \xi _{\ast
}\left( \boldsymbol{x}_{i}\right) >0\right\} .
\end{equation*}%
Then $c\left( \boldsymbol{\theta }\right) >0$ and $\ddot{k}\left( 
\boldsymbol{\theta }\right) \succeq c\left( \boldsymbol{\theta }\right)
\int_{\mathcal{\chi }}\boldsymbol{f}\left( \boldsymbol{x}\right) \boldsymbol{%
f}^{\T}\left( \boldsymbol{x}\right) \xi
_{\ast }\left( d\boldsymbol{x}\right) \succ \boldsymbol{0},$ so that $%
\boldsymbol{\theta }_{0}$ is the unique stationary point of the strictly
convex function $k\left( \boldsymbol{\theta }\right) $ and $\boldsymbol{0}_p$
is the unique zero of $\left\Vert \boldsymbol{\Psi }\left( \boldsymbol{%
\theta }\right) \right\Vert $.

Since $\boldsymbol{\Theta }$ is compact, its closed subset $\boldsymbol{%
\Theta }_{\delta }$ defined by $\left\Vert \boldsymbol{\theta }\right\Vert
\geq \delta $ is also compact, hence the continuous function $\left\Vert 
\boldsymbol{\Psi }\left( \boldsymbol{\theta }\right) \right\Vert $ attains
its $\inf $ there:%
\begin{equation*}
\inf_{\boldsymbol{\Theta }_{\delta }}\left\Vert \boldsymbol{\Psi }\left( 
\boldsymbol{\theta }\right) \right\Vert =\min_{\boldsymbol{\Theta }_{\delta
}}\left\Vert \boldsymbol{\Psi }\left( \boldsymbol{\theta }\right)
\right\Vert >0.
\end{equation*}%
This proves (ii), and completes the proof of Lemma \ref{lemma: limit}%
.\hfill\ $\square $

\medskip

\subsection{Proof of Theorem \protect\ref{thm:main}}

We must show that $\sqrt{n}\boldsymbol{M}_{n}^{1/2}\left( \boldsymbol{\hat{%
\theta}}_{n}-\boldsymbol{\tilde{\theta}}_{n}\right) \overset{d}{%
\longrightarrow }N\left( \boldsymbol{0}_{p},\sigma _{M}^{2}\boldsymbol{I}%
_{p}\right) $. The stochastic component of the estimate will be expressed in
terms of 
\begin{equation}
\boldsymbol{\bar{z}}_{n}\overset{def}{=}n^{-1}\sum_{i,j}\psi _{\sigma
}\left( \varepsilon _{ij}\right) \boldsymbol{f}\left( \boldsymbol{x}%
_{i}\right) .  \label{znbar}
\end{equation}
Note that $\boldsymbol{\bar{z}}_{n}$ is an average of independent r.v.s,
each with mean zero, and covariance matrices%
\begin{equation*}
\text{\textsc{var}}\left[ \psi _{\sigma }\left( \varepsilon \right) \right] 
\boldsymbol{f}\left( \boldsymbol{x}_{i}\right) \boldsymbol{f}^{\T}\left( \boldsymbol{x}_{i}\right) =E\left[ \psi
_{\sigma }^{2}\left( \varepsilon \right) \right] \boldsymbol{f}\left( 
\boldsymbol{x}_{i}\right) \boldsymbol{f}^{\T}\left( \boldsymbol{x}_{i}\right) ,
\end{equation*}%
so that \textsc{cov}$\left[ \boldsymbol{\bar{z}}_{n}\right] =n^{-1}E\left[
\psi _{\sigma }^{2}\left( \varepsilon \right) \right] \boldsymbol{M}%
_{n}\succ \boldsymbol{0}_{p}$. By the Lindberg-Feller CLT, 
\begin{equation}
\boldsymbol{M}_{n}^{-1/2}\sqrt{n}\boldsymbol{\bar{z}}_{n}\overset{d}{%
\rightarrow }N\left( \boldsymbol{0}_{p},E\left[ \psi _{\sigma }^{2}\left(
\varepsilon \right) \right] \boldsymbol{I}_{p}\right) .  \label{AN-z}
\end{equation}

Completing the proof of Theorem \ref{thm:main} requires a preliminary Lemma.

\begin{lemma}
\label{root-n consistent} Recall that $\boldsymbol{\tilde{\theta}}_{n}=%
\boldsymbol{\theta }_{0}+\boldsymbol{M}_{n}^{-1}\boldsymbol{b}_{n}$ . We
have that $\boldsymbol{\hat{\theta}}_{n}-\boldsymbol{\tilde{\theta}}_{n}$ is
a $\sqrt{n}$-consistent estimate of $\boldsymbol{0}_{p}$: 
\begin{equation}
\boldsymbol{\hat{\theta}}_{n}-\boldsymbol{\tilde{\theta}}_{n}=O_{p}\left(
n^{-1/2}\right) .  \label{claim2}
\end{equation}
\end{lemma}

\noindent \textit{Proof}: We use the functions 
\begin{eqnarray*}
a_{n}(\boldsymbol{x}) &=&\boldsymbol{f}^{\T}\left( \boldsymbol{x}\right) \left( \boldsymbol{\tilde{\theta}}_{n}-%
\boldsymbol{\theta }_{0}\right) =\boldsymbol{f}^{\T}\left( \boldsymbol{x}\right) \boldsymbol{M}_{n}^{-1}%
\boldsymbol{b}_{n}, \\
h_{n}\left( \boldsymbol{x}\right) &=&\tau \left( \boldsymbol{x}\right)
-a_{n}(\boldsymbol{x}), \\
d_{n}(\boldsymbol{x}) &=&\boldsymbol{f}^{\T}\left( \boldsymbol{x}\right) \left( \boldsymbol{\hat{\theta}}_{n}-%
\boldsymbol{\tilde{\theta}}_{n}\right) .
\end{eqnarray*}%
By (\ref{bound}) and C2), both $a_{n}(\boldsymbol{x})$ and $h_{n}\left( 
\boldsymbol{x}\right) $ are $o_{p}\left( 1\right) $.

We first show that $\boldsymbol{\hat{\theta}}_{n}-\boldsymbol{\tilde{\theta}}%
_{n}$ is a consistent estimate of $\boldsymbol{0}_{p}$: 
\begin{equation}
\boldsymbol{\hat{\theta}}_{n}-\boldsymbol{\tilde{\theta}}_{n}=o_{p}\left(
1\right) .  \label{claim1}
\end{equation}%
We have that $\Psi _{n}\left( \boldsymbol{\hat{\theta}}_{n}-\boldsymbol{%
\tilde{\theta}}_{n}\right) =n^{-1}\sum_{i,j}\psi _{\sigma }\left(
Y_{j}\left( \boldsymbol{x}_{i}\right) -\boldsymbol{f}^{\T}\left( \boldsymbol{x}_{i}\right) \boldsymbol{\hat{%
\theta}}_{n}\right) \boldsymbol{f}\left( \boldsymbol{x}_{i}\right) =%
\boldsymbol{0}_{p}$, by (\ref{m-eqn}). Now Theorem 5.9 of \cite{vw98}, with $%
\boldsymbol{\theta }_{0}=$ $\boldsymbol{0}_{p}$ and under the conditions
there -- ensured by our Lemma \ref{lemma: limit} -- asserts (\ref{claim1}),
and also that $d_{n}(\boldsymbol{x})=o_{p}\left( 1\right) $.

We require the stronger $\sqrt{n}$-consistency. For this, observe that by (%
\ref{m-eqn}), 
\begin{eqnarray}
\boldsymbol{0}_{p} &=&\frac{1}{n}\sum_{i,j}\psi _{\sigma }\left( Y_{j}\left( 
\boldsymbol{x}_{i}\right) -\boldsymbol{f}^{\T}\left( \boldsymbol{x}_{i}\right) \boldsymbol{\hat{\theta}}_{n}\right) 
\boldsymbol{f}\left( \boldsymbol{x}_{i}\right)  \notag \\
&=&\frac{1}{n}\sum_{i,j}\psi _{\sigma }\left( Y_{j}\left( \boldsymbol{x}%
_{i}\right) -\boldsymbol{f}^{\T}\left( 
\boldsymbol{x}_{i}\right) \boldsymbol{\tilde{\theta}}_{n}-d_{n}(\boldsymbol{x%
}_{i})\right) \boldsymbol{f}\left( \boldsymbol{x}_{i}\right)  \notag \\
&=&\frac{1}{n}\sum_{i,j}\psi _{\sigma }\left( Y_{j}\left( \boldsymbol{x}%
_{i}\right) -\boldsymbol{f}^{\T}\left( 
\boldsymbol{x}_{i}\right) \boldsymbol{\tilde{\theta}}_{n}\right) \boldsymbol{%
f}\left( \boldsymbol{x}_{i}\right) -\boldsymbol{G}_{n}\left( \boldsymbol{%
\hat{\theta}}_{n}-\boldsymbol{\tilde{\theta}}_{n}\right) ,  \label{exp1}
\end{eqnarray}%
where for generic $t_{i}\in \lbrack 0,1]$ and $u_{ij}\in \lbrack 0,1]$, and
with $s_{n,i}=h_{n}(\boldsymbol{x}_{i})-t_{i}d_{n}(\boldsymbol{x}%
_{i})=o_{p}\left( 1\right) $, 
\begin{eqnarray*}
\boldsymbol{G}_{n} &=&\frac{1}{n}\sum_{i,j}\psi _{\sigma }^{\prime }\left(
Y_{j}\left( \boldsymbol{x}_{i}\right) -\boldsymbol{f}^{\T}\left( \boldsymbol{x}_{i}\right) \boldsymbol{\tilde{%
\theta}}_{n}-t_{i}d_{n}(\boldsymbol{x}_{i})\right) \boldsymbol{f}\left( 
\boldsymbol{x}_{i}\right) \boldsymbol{f}^{\T}%
\left( \boldsymbol{x}_{i}\right) \\
 &=&\frac{1}{n}\sum_{i,j}\psi_{\sigma }^{\prime }
\left( \varepsilon_{ij}+s_{n,i}\right) \boldsymbol{f}\left( \boldsymbol{x}_{i}\right) 
\boldsymbol{f}^{\T}\left( \boldsymbol{x}%
_{i}\right) \\
&=&\frac{1}{n}\sum_{i,j}\psi _{\sigma }^{\prime }\left( \varepsilon
_{ij}\right) \boldsymbol{f}\left( \boldsymbol{x}_{i}\right) \boldsymbol{f}^{\T}\left( \boldsymbol{x}_{i}\right) +\frac{1}{n}%
\sum_{i,j}\psi _{\sigma }^{\prime \prime }\left( \varepsilon
_{ij}+u_{ij}s_{n,i}\right) s_{n,i}\boldsymbol{f}\left( \boldsymbol{x}%
_{i}\right) \boldsymbol{f}^{\T}\left( 
\boldsymbol{x}_{i}\right) .
\end{eqnarray*}%
We will show that $\boldsymbol{G}_{n}$ tends in probability to the positive
definite matrix $\boldsymbol{G}_{\ast }=E\left[ \psi _{\sigma }^{\prime
}\left( \varepsilon \right) \right] \int_{\mathcal{\chi }}\boldsymbol{f}%
\left( \boldsymbol{x}\right) \boldsymbol{f}^{\T}
\left( \boldsymbol{x}\right) \xi _{\ast }\left( d\boldsymbol{x}\right) $.
Then by (\ref{exp1}), and ignoring terms which are $o_p(n^{-1/2}),$%
\begin{eqnarray*}
\sqrt{n}\left( \boldsymbol{\hat{\theta}}_{n}-\boldsymbol{\tilde{\theta}}%
_{n}\right) &=&\left( \boldsymbol{G}_{\ast }^{-1}+o_{p}(1)\right) \left\{ 
\frac{1}{\sqrt{n}}\sum_{i,j}\psi _{\sigma }\left( Y_{j}\left( \boldsymbol{x}%
_{i}\right) -\boldsymbol{f}^{\T}\left( 
\boldsymbol{x}_{i}\right) \boldsymbol{\tilde{\theta}}_{n}\right) \boldsymbol{%
f}\left( \boldsymbol{x}_{i}\right) \right\} \\
&=&\left( \boldsymbol{G}_{\ast }^{-1}+o_{p}(1)\right) \left\{ \frac{1}{\sqrt{%
n}}\sum_{i,j}\psi _{\sigma }\left( \varepsilon _{ij}+h_{n}\left( \boldsymbol{%
x}_{i}\right) \right) \boldsymbol{f}\left( \boldsymbol{x}_{i}\right) \right\}
\\
&=&\left( \boldsymbol{G}_{\ast }^{-1}+o_{p}(1)\right) \left\{ \sqrt{n}%
\boldsymbol{\bar{z}}_{n}+\frac{1}{n}\sum_{i,j}\psi _{\sigma }^{\prime
}\left( \varepsilon _{ij}+t_{i}h_{n}\left( \boldsymbol{x}_{i}\right) \right) 
\sqrt{n}h\left( \boldsymbol{x}_{i}\right) \boldsymbol{f}\left( \boldsymbol{x}%
_{i}\right) \right\} ,
\end{eqnarray*}%
where the first term in the final set of braces is $O_{p}(1)$ by (\ref{AN-z}%
), and the second is $O_{p}(1)$ because 
\begin{equation*}
\sqrt{n}h_{n}\left( \boldsymbol{x}\right) =\sqrt{n}\tau \left( \boldsymbol{x}%
\right) -\sqrt{n}\boldsymbol{f}^{\T}\left( 
\boldsymbol{x}\right) \boldsymbol{M}_{n}^{-1}\boldsymbol{b}_{n}=\sqrt{n}\tau
\left( \boldsymbol{x}\right) -\boldsymbol{f}^{\T}\left( \boldsymbol{x}\right) \boldsymbol{M}_{n}^{-1}\int_{\mathcal{\chi }}%
\boldsymbol{f}\left( \boldsymbol{x}\right) \sqrt{n}\tau \left( \boldsymbol{x}%
\right) \xi _{n}\left( d\boldsymbol{x}\right)
\end{equation*}%
is bounded, again by virtue of (\ref{bound}). Thus showing that 
\begin{equation}
\boldsymbol{G}_{n}\overset{pr}{\rightarrow }\boldsymbol{G}_{\ast },
\label{A-limit}
\end{equation}%
will complete the verification of (\ref{claim2}).

To show (\ref{A-limit}) note that by the triangle inequality, 
\begin{equation}
\left\Vert \boldsymbol{G}_{n}-\boldsymbol{G}_{\ast }\right\Vert \leq
\left\Vert \boldsymbol{G}_{n}-\frac{1}{n}\sum_{i,j}\psi _{\sigma }^{\prime
}\left( \varepsilon _{ij}\right) \boldsymbol{f}\left( \boldsymbol{x}%
_{i}\right) \boldsymbol{f}^{\T}\left( 
\boldsymbol{x}_{i}\right) \right\Vert +\left\Vert \boldsymbol{W}%
_{n}\right\Vert +\left\Vert \frac{1}{n}\sum_{i,j}E\left[ \psi _{\sigma
}^{\prime }\left( \varepsilon _{ij}\right) \right] \boldsymbol{f}\left( 
\boldsymbol{x}_{i}\right) \boldsymbol{f}^{\T}%
\left( \boldsymbol{x}_{i}\right) -\boldsymbol{G}_{\ast }\right\Vert ,
\label{A-decomp}
\end{equation}%
where 
\begin{equation*}
\boldsymbol{W}_{n}\overset{def}{=}\frac{1}{n}\sum_{i,j}\left( \psi _{\sigma
}^{\prime }\left( \varepsilon _{ij}\right) -E\left[ \psi _{\sigma }^{\prime
}\left( \varepsilon _{ij}\right) \right] \right) \boldsymbol{f}\left( 
\boldsymbol{x}_{i}\right) \boldsymbol{f}^{\T}\left( \boldsymbol{x}_{i}\right) =\sum_{i}\xi _{n}\left( \boldsymbol{x}%
_{i}\right) V_{i}\boldsymbol{f}\left( \boldsymbol{x}_{i}\right) \boldsymbol{f%
}^{\T}\left( \boldsymbol{x}_{i}\right) ,
\end{equation*}%
for independent r.v.s $V_{i}=\sum_{j=1}^{n_{i}}\left( \psi _{\sigma
}^{\prime }\left( \varepsilon _{ij}\right) -E\left[ \psi _{\sigma }^{\prime
}\left( \varepsilon _{ij}\right) \right] \right) /n_{i}$ with means $0$ and
variances \textsc{var}$\left[ \psi _{\sigma }^{\prime }\left( \varepsilon
\right) \right] /n_{i}\leq m_{1}^{2}/n_{i}$.

The first norm on the rhs of (\ref{A-decomp}) is 
\begin{equation*}
\left\Vert \frac{1}{n}\sum_{i,j}\psi _{\sigma }^{\prime \prime }\left(
\varepsilon _{ij}+u_{ij}s_{n,i}\right) s_{n,i}\boldsymbol{f}\left( 
\boldsymbol{x}_{i}\right) \boldsymbol{f}^{\T}\left( \boldsymbol{x}_{i}\right) \right\Vert \leq \frac{m_{2}}{n}%
\sum_{i,j}\left\vert s_{n,i}\right\vert \left\Vert \boldsymbol{f}\left( 
\boldsymbol{x}_{i}\right) \right\Vert ^{2}=o_{p}\left( 1\right) ,
\end{equation*}%
and the third is $\left\Vert \int_{\mathcal{\chi }}E\left[ \psi _{\sigma
}^{\prime }\left( \varepsilon \right) \right] \boldsymbol{f}\left( 
\boldsymbol{x}\right) \boldsymbol{f}^{\T}\left( 
\boldsymbol{x}\right) \left( \xi _{n}-\xi _{\ast }\right) \left( d%
\boldsymbol{x}\right) \right\Vert =o(1)$. A typical element of $W_{n}$ is $%
W_{n,k,l}=\sum_{i}\xi _{n}\left( \boldsymbol{x}_{i}\right) V_{i}f_{k}\left( 
\boldsymbol{x}_{i}\right) f_{l}\left( \boldsymbol{x}_{i}\right) $, with mean
zero and 
\begin{equation*}
\text{\textsc{var}}\left[ W_{n,k,l}\right] =\sum_{i}\frac{n_{i}^{2}}{n^{2}}%
\text{\textsc{var}}\left[ V_{i}\right] \left( f_{k}\left( \boldsymbol{x}%
_{i}\right) f_{l}\left( \boldsymbol{x}_{i}\right) \right) ^{2}\leq \frac{%
m_{1}^{2}}{n}\int_{\mathcal{\chi }}\left( f_{k}\left( \boldsymbol{x}%
_{i}\right) f_{l}\left( \boldsymbol{x}_{i}\right) \right) ^{2}\xi _{n}\left(
d\boldsymbol{x}\right) =O(n^{-1}),
\end{equation*}%
so that by Chebyshev's Inequality, $W_{n,k,l}\overset{pr}{\rightarrow }0$
and so $\boldsymbol{W}_{n}\overset{pr}{\rightarrow }\boldsymbol{0}.$ Now (%
\ref{A-limit}) follows, completing the proof of the Lemma.\hfill\ $\square
\medskip $

We can now complete the proof of Theorem \ref{thm:main}. To establish the
asymptotic normality of $\boldsymbol{\hat{\theta}}_{n}-\boldsymbol{\tilde{%
\theta}}_{n}$, we observe that by (\ref{m-eqn}), 
\begin{eqnarray}
\boldsymbol{0}_p &=& \frac{1}{n}\sum_{i,j}\psi _{\sigma }\left( Y_{j}\left( 
\boldsymbol{x}_{i}\right) -\boldsymbol{f}^{\T}\left( \boldsymbol{x}_{i}\right) \boldsymbol{\hat{\theta}}_{n}\right) 
\boldsymbol{f}\left( \boldsymbol{x}_{i}\right)  \notag \\
&=&\frac{1}{n}\sum_{i,j}\psi _{\sigma }\left( Y_{j}\left( \boldsymbol{x}%
_{i}\right) -\boldsymbol{f}^{\T}\left( 
\boldsymbol{x}_{i}\right) \boldsymbol{\tilde{\theta}}_{n}-d_{n}\left( 
\boldsymbol{x}_{i}\right) \right) \boldsymbol{f}\left( \boldsymbol{x}%
_{i}\right)  \notag \\
\text{ }& =&\frac{1}{n}\sum_{i,j}\psi _{\sigma }\left( Y_{j}\left( 
\boldsymbol{x}_{i}\right) -\boldsymbol{f}^{\T}\left( \boldsymbol{x}_{i}\right) \boldsymbol{\tilde{\theta}}_{n}\right) 
\boldsymbol{f}\left( \boldsymbol{x}_{i}\right)  \notag \\
\text{}&&-\frac{1}{n}\sum_{i,j}\psi _{\sigma }^{\prime }\left( Y_{j}\left( 
\boldsymbol{x}_{i}\right) -\boldsymbol{f}^{\T}\left( \boldsymbol{x}_{i}\right) \boldsymbol{\tilde{\theta}}_{n}\right) 
\boldsymbol{f}\left( \boldsymbol{x}_{i}\right) d_{n}\left( \boldsymbol{x}%
_{i}\right) +\boldsymbol{S}_{n},  \label{exp2}
\end{eqnarray}%
where 
\begin{equation*}
\boldsymbol{S}_{n}=\left( 1/2n\right) \sum_{i,j}\psi _{\sigma }^{\prime
\prime }\left( Y_{j}\left( \boldsymbol{x}_{i}\right) -\boldsymbol{f}^{\T}\left( \boldsymbol{x}_{i}\right) \boldsymbol{%
\tilde{\theta}}_{n}+t_id_{n}\left( \boldsymbol{x}_{i}\right) \right)
d_{n}^{2}\left( \boldsymbol{x}_{i}\right) \boldsymbol{f}\left( \boldsymbol{x}%
_{i}\right) .
\end{equation*}%
With $B_{n}\overset{def}{=}ch_{\max }\int_{\mathcal{\chi }}\boldsymbol{f}%
\left( \boldsymbol{x}\right) \boldsymbol{f}^{\T}%
\left( \boldsymbol{x}\right) \left\Vert \boldsymbol{f}\left( \boldsymbol{x}%
\right) \right\Vert \xi _{n}\left( d\boldsymbol{x}\right) $ $\leq $ $\int_{%
\mathcal{\chi }}\left\Vert \boldsymbol{f}\left( \boldsymbol{x}\right)
\right\Vert ^{3}\xi _{n}\left( d\boldsymbol{x}\right) <\infty $, we have 
\begin{equation*}
\left\Vert \boldsymbol{S}_{n}\right\Vert \leq \frac{m_{2}}{2n}%
\sum_{i,j}\left( \boldsymbol{\hat{\theta}}_{n}-\boldsymbol{\tilde{\theta}}%
_{n}\right) ^{\T}\boldsymbol{f}\left( 
\boldsymbol{x}_{i}\right) \boldsymbol{f}^{\T}\left( \boldsymbol{x}_{i}\right) \left( \boldsymbol{\hat{\theta}}_{n}-%
\boldsymbol{\tilde{\theta}}_{n}\right) \left\Vert \boldsymbol{f}\left( 
\boldsymbol{x}_{i}\right) \right\Vert \leq \frac{m_{2}}{2}B_{n}\left\Vert 
\boldsymbol{\hat{\theta}}_{n}-\boldsymbol{\tilde{\theta}}_{n}\right\Vert
^{2}=O_{p}\left( n^{-1}\right) .
\end{equation*}%
The final step in the preceding line uses Lemma \ref{root-n consistent}.
Thus, from (\ref{exp2}) and disregarding terms that are $O_{p}\left(
n^{-1/2}\right) $, 
\begin{equation}
\sqrt{n}\left( \boldsymbol{\hat{\theta}}_{n}-\boldsymbol{\tilde{\theta}}%
_{n}\right) =\boldsymbol{P}_{n}^{-1}\frac{1}{\sqrt{n}}\sum_{i,j}\psi
_{\sigma }\left( Y_{j}\left( \boldsymbol{x}_{i}\right) -\boldsymbol{f}^{\T}\left( \boldsymbol{x}_{i}\right) \boldsymbol{%
\tilde{\theta}}_{n}\right) \boldsymbol{f}\left( \boldsymbol{x}_{i}\right) ,
\label{exp3}
\end{equation}%
where%
\begin{eqnarray}
\boldsymbol{P}_{n} &\overset{def}{=}&\frac{1}{n}\sum_{i,j}\psi _{\sigma
}^{\prime }\left( Y_{j}\left( \boldsymbol{x}_{i}\right) -\boldsymbol{f}^{\T}\left( \boldsymbol{x}_{i}\right) \boldsymbol{%
\tilde{\theta}}_{n}\right) \boldsymbol{f}\left( \boldsymbol{x}_{i}\right) 
\boldsymbol{f}^{\T}\left( \boldsymbol{x}%
_{i}\right)  \notag \\
&=&\frac{1}{n}\sum_{i,j}\psi _{\sigma }^{\prime }\left( \tau \left( 
\boldsymbol{x}_{i}\right) +\varepsilon _{ij}-a_{n}\left( \boldsymbol{x}%
_{i}\right) \right) \boldsymbol{f}\left( \boldsymbol{x}_{i}\right) 
\boldsymbol{f}^{\T}\left( \boldsymbol{x}%
_{i}\right)  \notag \\
&=&\frac{1}{n}\sum_{i,j}\psi _{\sigma }^{\prime }\left( \tau \left( 
\boldsymbol{x}_{i}\right) +\varepsilon _{ij}\right) \boldsymbol{f}\left( 
\boldsymbol{x}_{i}\right) \boldsymbol{f}^{\T}\left( \boldsymbol{x}_{i}\right) -\boldsymbol{T}_{n},  \label{exp4}
\end{eqnarray}%
and 
\begin{equation*}
\boldsymbol{T}_{n}=\frac{1}{n}\sum_{i,j}\psi _{\sigma }^{\prime \prime
}\left( \tau \left( \boldsymbol{x}_{i}\right) +\varepsilon
_{ij}-t_ia_{n}\left( \boldsymbol{x}_{i}\right) \right) a_{n}\left( 
\boldsymbol{x}_{i}\right) \boldsymbol{f}\left( \boldsymbol{x}_{i}\right) 
\boldsymbol{f}^{\T}\left( \boldsymbol{x}%
_{i}\right) .
\end{equation*}%
For any vector $\boldsymbol{\alpha }$ of unit norm, 
\begin{eqnarray*}
\left\vert \boldsymbol{\alpha }^{\T}\boldsymbol{%
T}_{n}\boldsymbol{\alpha }\right\vert &\leq &\frac{1}{n}\sum_{i,j}\left\vert
\psi _{\sigma }^{\prime \prime }\left( \tau \left( \boldsymbol{x}_{i}\right)
+\varepsilon _{ij}-ta_{n}\left( \boldsymbol{x}_{i}\right) \right)
a_{n}\left( \boldsymbol{x}_{i}\right) \right\vert \left( \boldsymbol{\alpha }%
^{\T}\boldsymbol{f}\left( \boldsymbol{x}%
_{i}\right) \right) ^{2} \\
&\leq &\frac{m_{2}}{n}\sum_{i,j}\left( \boldsymbol{\alpha }^{\T}\boldsymbol{f}\left( \boldsymbol{x}_{i}\right)
\right) ^{2}\left( \max_{\mathcal{\chi }}\left\vert a_{n}\left( \boldsymbol{x%
}\right) \right\vert \right) \\
&=&m_{2}\cdot \boldsymbol{\alpha }^{\T}%
\boldsymbol{M}_{n}\boldsymbol{\alpha }\cdot \max_{\mathcal{\chi }}\left\Vert 
\boldsymbol{f}\left( \boldsymbol{x}\right) \right\Vert \cdot \left\Vert 
\boldsymbol{M}_{n}^{-1}\boldsymbol{b}_{n}\right\Vert ,
\end{eqnarray*}%
which is $O\left( n^{-1/2}\right) $. Thus all eigenvalues of $\boldsymbol{T}%
_{n}$ are $O\left( n^{-1/2}\right) $, hence so is $\boldsymbol{T}_{n}$
itself and (\ref{exp4}) becomes 
\begin{equation*}
\boldsymbol{P}_{n}=E\left[ \psi _{\sigma }^{\prime }\left( \varepsilon
\right) \right] \boldsymbol{M}_{n}-O_p\left( n^{-1/2}\right) ,
\end{equation*}%
where we invoke the WLLN. Now (\ref{exp3}) becomes%
\begin{equation}
\sqrt{n}\left( \boldsymbol{\hat{\theta}}_{n}-\boldsymbol{\tilde{\theta}}%
_{n}\right) =\left[ E\left[ \psi _{\sigma }^{\prime }\left( \varepsilon
\right) \right] \boldsymbol{M}_{n}+O\left( n^{-1/2}\right) \right]
^{-1}\cdot \sqrt{n}\boldsymbol{\tilde{z}}_{n},  \label{exp5}
\end{equation}%
for%
\begin{equation*}
\boldsymbol{\tilde{z}}_{n}=\frac{1}{n}\sum_{i,j}\psi _{\sigma }\left(
Y_{j}\left( \boldsymbol{x}_{i}\right) -\boldsymbol{f}^{\T}\left( \boldsymbol{x}_{i}\right) \boldsymbol{\tilde{%
\theta}}_{n}\right) \boldsymbol{f}\left( \boldsymbol{x}_{i}\right) =\frac{1}{%
n}\sum_{i,j}\psi _{\sigma }\left( \varepsilon _{ij}+h_{n}\left( \boldsymbol{x%
}_{i}\right) \right) \boldsymbol{f}\left( \boldsymbol{x}_{i}\right) .
\end{equation*}%
Recalling that $h_{n}\left( \boldsymbol{x}\right) $ is $O\left(
n^{-1/2}\right) $, and using (\ref{znbar}), we have%
\begin{equation*}
\boldsymbol{\tilde{z}}_{n}=\frac{1}{n}\sum_{i,j}\psi _{\sigma }\left(
\varepsilon _{ij}+h_{n}\left( \boldsymbol{x}_{i}\right) \right) \boldsymbol{f%
}\left( \boldsymbol{x}_{i}\right) =\boldsymbol{\bar{z}}_{n}+\frac{1}{n}%
\sum_{i,j}\psi _{\sigma }^{\prime }\left( \varepsilon _{ij}\right)
h_{n}\left( \boldsymbol{x}_{i}\right) \boldsymbol{f}\left( \boldsymbol{x}%
_{i}\right) +O_{p}\left( n^{-1}\right) ,
\end{equation*}%
with (and again invoking the WLLN)%
\begin{eqnarray*}
\sqrt{n}\boldsymbol{\tilde{z}}_{n}-\sqrt{n}\boldsymbol{\bar{z}}_{n} &=&\frac{%
1}{n}\sum_{i,j}\psi _{\sigma }^{\prime }\left( \varepsilon _{ij}\right) 
\sqrt{n}h_{n}\left( \boldsymbol{x}_{i}\right) \boldsymbol{f}\left( 
\boldsymbol{x}_{i}\right) +O_{p}\left( n^{-1/2}\right) \\
&=&E\left[ \psi _{\sigma }^{\prime }\left( \varepsilon \right) \right] \int_{%
\mathcal{\chi }}\sqrt{n}h_{n}\left( \boldsymbol{x}\right) \boldsymbol{f}%
\left( \boldsymbol{x}\right) \xi _{n}\left( d\boldsymbol{x}\right)
+o_{p}\left( 1\right) \\
&=&E\left[ \psi _{\sigma }^{\prime }\left( \varepsilon \right) \right]
\left\{ \sqrt{n}\boldsymbol{b}_{n}-\sqrt{n}\boldsymbol{M}_{n}\left( 
\boldsymbol{M}_{n}^{-1}\boldsymbol{b}_{n}\right) \right\} +o_{p}\left(
1\right) ,
\end{eqnarray*}%
which is $o_{p}\left( 1\right) $. Thus $\sqrt{n}\boldsymbol{\tilde{z}}_{n}=%
\sqrt{n}\boldsymbol{\bar{z}}_{n}+o_{p}(1)$ and (\ref{exp5}) becomes%
\begin{equation*}
\sqrt{n}\left( \boldsymbol{\hat{\theta}}_{n}-\boldsymbol{\tilde{\theta}}%
_{n}\right) =\left[ E\left[ \psi _{\sigma }^{\prime }\left( \varepsilon
\right) \right] \boldsymbol{M}_{n}+O\left( n^{-1/2}\right) \right]
^{-1}\left( \sqrt{n}\boldsymbol{\bar{z}}_{n}+o_{p}\left( 1\right) \right) ,
\end{equation*}%
whence 
\begin{equation}
\sqrt{n}\boldsymbol{M}_{n}^{1/2}\left( \boldsymbol{\hat{\theta}}_{n}-%
\boldsymbol{\tilde{\theta}}_{n}\right) =\frac{1}{E\left[ \psi _{\sigma }^{{%
\prime} }\left( \varepsilon \right) \right] }\boldsymbol{M}_{n}^{-1/2}\sqrt{n%
}\boldsymbol{\bar{z}}_{n}+o_{p}\left( 1\right) .  \label{represent}
\end{equation}%
From this, the statement of the Theorem is immediate.\hfill\ $\square $

\subsection{Proof of Theorem \protect\ref{thm: max loss}}

To maximize $\textsc{aimse}$ over $\tau \in \Upsilon $, note that both of
its first two terms at (\ref{IMSEn}) become larger if $\tau \left( 
\boldsymbol{x}\right) $ is multiplied by a constant exceeding one in
absolute value, hence at a maximum (\ref{bound}) is attained with equality.
Define $\tau _{0}\left( \boldsymbol{x}\right) =\sqrt{n}\tau \left( 
\boldsymbol{x}\right) /\kappa \text{ and } \boldsymbol{c} \left( \xi \right)
= \sqrt{n} \boldsymbol{b}_0 \left( \xi \right) / \kappa. $ Then with $%
\boldsymbol{c}_{n}=\boldsymbol{c}\left( \xi _{n}\right) $ we have that $%
\max_{\tau \in \Upsilon }\textsc{aimse}\left( \xi _{n}|\tau \right) $ is 
\begin{equation}
J\left( \xi _{n}\right) \overset{def}{=}\sigma _{M}^{2}tr\boldsymbol{AM}%
_{n}^{-1}+\kappa ^{2}\left( 1+\max_{\tau _{0}}\boldsymbol{c}_{n}^{\T}\boldsymbol{M}_{n}^{-1}\mathbf{A}\boldsymbol{M}%
_{n}^{-1}\boldsymbol{c}_{n}\right) ,  \label{I1}
\end{equation}%
with $\tau _{0}$ -- which appears in $\boldsymbol{c}_{n}$ via $\boldsymbol{b}%
_{0}$ -- constrained by (\ref{orthogonality}) and by $\sum_{i=1}^{N}\tau
_{0}^{2}\left( \boldsymbol{x}_{i}\right) =1$.

To carry out the maximization in (\ref{I1}) recall that a feature of the
Gram-Schmidt process used in the construction of $\boldsymbol{Q}$ is that $%
\boldsymbol{F}=\boldsymbol{QT}$ for a non-singular lower triangular matrix $%
\boldsymbol{T}$. Let $\boldsymbol{Q}_{\boldsymbol{\perp }}:N\times N-p$ be
the orthogonal complement of $\boldsymbol{Q}$, so that $\left( \boldsymbol{Q}%
\vdots \boldsymbol{Q}_{\boldsymbol{\perp }}\right) :N\times N$ is
orthogonal. The condition in (\ref{orthogonality}) asserts that $\boldsymbol{%
\tau }_{0}$ is orthogonal to the columns of $\boldsymbol{Q}$, hence is a
linear combination of the columns of $\boldsymbol{Q}_{\boldsymbol{\perp }}$
and so there is $\boldsymbol{\beta }_{N-p\times 1}$ of unit norm for which $%
\boldsymbol{\tau }_{0}=\boldsymbol{Q}_{\boldsymbol{\perp }}\boldsymbol{\beta 
}$.

In this notation $\boldsymbol{A}=\boldsymbol{F}^{\T}\boldsymbol{F}=\boldsymbol{T}^{\T}\boldsymbol{T}$ and 
\begin{eqnarray*}
\boldsymbol{M}_{n} &=&\boldsymbol{f}^{\T}%
\boldsymbol{D}\left( \xi _{n}\right) \boldsymbol{F}=\boldsymbol{T}^{\T}\boldsymbol{Q}^{\T}\boldsymbol{D}\left( \xi _{n}\right) \boldsymbol{QT}%
, \\
\boldsymbol{c}_{n} &=&\boldsymbol{f}^{\T}%
\boldsymbol{D}\left( \xi _{n}\right) \boldsymbol{\tau }_{0}=\boldsymbol{T}^{\T}\boldsymbol{Q}^{\T}\boldsymbol{D}\left( \xi _{n}\right) \boldsymbol{Q}%
_{\boldsymbol{\perp }}{\boldsymbol{\beta }};
\end{eqnarray*}%
then (\ref{I1}) becomes%
\begin{equation}
J\left( \xi _{n}\right) =\sigma _{M}^{2}tr\left( \boldsymbol{Q}^{\T}\boldsymbol{D}\left( \xi _{n}\right) \boldsymbol{Q}%
\right) ^{-1}+\kappa ^{2}\left( 1+\max_{\left\Vert \boldsymbol{\beta }%
\right\Vert =1}\boldsymbol{\beta }^{\T}%
\boldsymbol{Q}_{\boldsymbol{\perp }}^{\T}%
\boldsymbol{D}\left( \xi _{n}\right) \boldsymbol{Q}\left( \boldsymbol{Q}^{\T}\boldsymbol{D}\left( \xi _{n}\right) 
\boldsymbol{Q}\right) ^{-2}\boldsymbol{Q}^{\T}\boldsymbol{D}\left( \xi _{n}\right) \boldsymbol{Q}_{\boldsymbol{\perp }}{%
\boldsymbol{\beta }}\right) .  \label{I2}
\end{equation}%
Noting that the maximum eigenvalue of a matrix $\boldsymbol{PP}^{\T}$ is that of $\boldsymbol{P}^{\T}\boldsymbol{P}$, and that $\boldsymbol{Q}_{%
\boldsymbol{\perp }}\boldsymbol{Q}_{\boldsymbol{\perp }}^{\T}=\left( \boldsymbol{I}_{N}-\boldsymbol{QQ}^{\T}\right) $, we have 
\begin{eqnarray*}
&&\max_{\left\Vert \boldsymbol{\beta }\right\Vert =1}\boldsymbol{\beta }^{\T}\boldsymbol{Q}_{\boldsymbol{\perp }}^{\T}\boldsymbol{D}\left( \xi _{n}\right) 
\boldsymbol{Q}\left( \boldsymbol{Q}^{\T}%
\boldsymbol{D}\left( \xi _{n}\right) \boldsymbol{Q}\right) ^{-2}\boldsymbol{Q%
}^{\T}\boldsymbol{D}\left( \xi _{n}\right) 
\boldsymbol{Q}_{\boldsymbol{\perp }}\boldsymbol{\beta } \\
&=&ch_{\max }\boldsymbol{Q}_{\boldsymbol{\perp }}^{\T}\boldsymbol{D}\left( \xi _{n}\right) \boldsymbol{Q}%
\left( \boldsymbol{Q}^{\T}\boldsymbol{D}%
\left( \xi _{n}\right) \boldsymbol{Q}\right) ^{-1}\cdot \left( \boldsymbol{Q}%
^{\T}\boldsymbol{D}\left( \xi _{n}\right) 
\boldsymbol{Q}\right) ^{-1}\boldsymbol{Q}^{\T}
\boldsymbol{D}\left( \xi _{n}\right) \boldsymbol{Q}_{\boldsymbol{\perp }} \\
&=&ch_{\max }\left( \boldsymbol{Q}^{\T}%
\boldsymbol{D}\left( \xi _{n}\right) \boldsymbol{Q}\right) ^{-1}\boldsymbol{Q%
}^{\T}\boldsymbol{D}\left( \xi _{n}\right) 
\boldsymbol{Q}_{\boldsymbol{\perp }}\boldsymbol{Q}_{\boldsymbol{\perp }}^{\T}\boldsymbol{D}\left( \xi _{n}\right) 
\boldsymbol{Q}\left( \boldsymbol{Q}^{\T}%
\boldsymbol{D}\left( \xi _{n}\right) \boldsymbol{Q}\right) ^{-1} \\
&=&ch_{\max }\left( \boldsymbol{Q}^{\T}%
\boldsymbol{D}\left( \xi _{n}\right) \boldsymbol{Q}\right) ^{-1}\boldsymbol{Q%
}^{\T}\boldsymbol{D}^{2}\left( \xi
_{n}\right) \boldsymbol{Q}\left( \boldsymbol{Q}^{\T}\boldsymbol{D}\left( \xi _{n}\right) \boldsymbol{Q}%
\right) ^{-1}-1.
\end{eqnarray*}%
This in (\ref{I2}) gives (\ref{I3}).\hfill\ $\square $

\subsection{Proof of Lemma \protect\ref{lemma: max IMloss} }

Using (\ref{AN}), an asymptotic approximation is 
\begin{equation*}
\text{\textsc{aim}}=\sum_{i=1}^{N}E\left[ \rho \left( \tau \left( \boldsymbol{%
x}_{i}\right) -\boldsymbol{f}^{\T}\left( 
\boldsymbol{x}_{i}\right) \left( \boldsymbol{M}_{n}^{-1}\boldsymbol{b}%
_{n}+\left( \sigma _{M}/\sqrt{n}\right) \boldsymbol{M}_{n}^{-1/2}\boldsymbol{%
Z}\right) \right) \right],
\end{equation*}%
for a random vector $\boldsymbol{Z}\sim \boldsymbol{N}\left( \boldsymbol{0}%
_{p},\boldsymbol{I}_{p}\right) $. Define random variables $W_{i}$ by \begin{equation*}
W_{i}=\boldsymbol{f}^{\T}\left( \boldsymbol{x}_{i}\right) \left( \sigma _{M}/\sqrt{n}\right) \boldsymbol{M}_{n}^{-1/2}%
\boldsymbol{Z}\sim N\left( 0,\left( \sigma _{M}^{2}/n\right) \boldsymbol{f}^{\T}\left( \boldsymbol{x}_{i}\right) \boldsymbol{M}_{n}^{-1}\boldsymbol{f}\left( \boldsymbol{x}_{i}\right) \right) , 
\end{equation*}
and set 
$\alpha _{i}\left( \boldsymbol{\tau }\right) =\tau \left( \boldsymbol{x}_{i}\right) -\boldsymbol{f}^{\T}\left( 
\boldsymbol{x}_{i}\right) \boldsymbol{M}_{n}^{-1}\boldsymbol{b}_{n}$. Then 
\begin{equation}
\text{\textsc{aim}}=\sum_{i=1}^{N}E\left[ \rho \left( \alpha _{i}\left( 
\boldsymbol{\tau }\right) -W_{i}\right) \right] .  \label{ghapp}
\end{equation}
To write \textsc{aim} as a function of $\boldsymbol{\beta }$ we note that, in terms of the orthogonal basis,
\begin{equation*}    
\boldsymbol{f}^{\T}\left( \boldsymbol{x}%
_{i}\right) \boldsymbol{M}_{n}^{-1}\boldsymbol{b}_{n}=\left( \kappa/\sqrt{n}\right)\left( \boldsymbol{Q\boldsymbol{R}}^{-1}\left( \boldsymbol{%
\xi }_{n}\right) \boldsymbol{Q}^{\T}\boldsymbol{D}\left( 
\boldsymbol{\xi }_{n}\right) \boldsymbol{Q}_{\perp }\boldsymbol{\beta }%
\right) _{i},
\end{equation*}
whence 
$
\alpha _{i}\left( \boldsymbol{\tau }\right) =\left( \kappa/\sqrt{n}\right)%
\boldsymbol{h}_{i}^{\T}\left( \boldsymbol{\xi }_{n}\right) \boldsymbol{\beta }$.
Since 
\begin{equation*}
\text{\textsc{var}}\left[ W_{i}\right] =\left( \sigma _{M}^{2}/n\right) 
\boldsymbol{f}^{\T}\left( \boldsymbol{x}%
_{i}\right) \boldsymbol{M}_{n}^{-1}\boldsymbol{f}\left( \boldsymbol{x}%
_{i}\right) =\left( \sigma _{M}^{2}/n\right) \left[ 
\boldsymbol{Q\boldsymbol{R}}^{-1}\left( \boldsymbol{\xi }_{n}\right) 
\boldsymbol{Q}^{\T}\right] _{ii},
\end{equation*}%
we can write $W_{i}=\left( \sigma _{M}/\sqrt{n}\right) s_i\left( \boldsymbol{\xi }_{n}\right)Z_{i}$.
This in (\ref{ghapp}) gives (\ref{im}). $\hfill
\square $

\subsection{Proof of Lemma \protect\ref{lemma: derivs} }

That \textsc{aim} is given by (\ref{G}) is immediate from (\ref{im}). Then
the $j^{th}$ element of $\boldsymbol{G}_{1}\left( \boldsymbol{\beta }\right) 
$ is 
\begin{eqnarray*}
\frac{\partial G\left( \boldsymbol{\beta }\right) }{\partial \beta _{j}}
&=&\sum_{i=1}^{N}\int_{-\infty }^{\infty }\psi \left( \frac{\kappa }{\sqrt{n}%
}\boldsymbol{h}_{i}^{\T}\left( \boldsymbol{\xi }%
_{n}\right) \boldsymbol{\beta }-\frac{\sigma _{M}}{\sqrt{n}}s_{i}\left( 
\boldsymbol{\xi }_{n}\right) z\right) \phi \left( z\right) dz\frac{\partial 
\boldsymbol{h}_{i}^{\T}\left( \boldsymbol{\xi }%
_{n}\right) \boldsymbol{\beta }}{\partial \beta _{j}} \\
&=&\frac{\kappa }{\sqrt{n}}\sum_{i=1}^{N}\int_{-\infty }^{\infty }
\left({\psi }\left(\boldsymbol{\eta}\left( \boldsymbol{\xi }_{n},z\right)\right)\right)_i 
\phi\left(z\right)dz H_{ij} \\
&=& \frac{\kappa }{\sqrt{n}}\left(
\boldsymbol{H}^{\T}\left( \boldsymbol{\xi }_{n}\right) 
\int_{-\infty }^{\infty }
{\psi}\left(\boldsymbol{\eta}\left( \boldsymbol{\xi }_{n},z\right) \right)
\phi \left( z\right) dz \right)_j,
\end{eqnarray*}
as claimed. The $\left( j,k\right) ^{th}$ element of $\boldsymbol{G}%
_{2}\left( \boldsymbol{\beta }\right) $ is%
\begin{eqnarray*}
\frac{\partial \left( \boldsymbol{G}_{1}\left( \boldsymbol{\beta }\right)
\right) _{j}}{\partial \beta _{k}} 
&=&\frac{\partial }{\partial \beta _k}%
\frac{\kappa }{\sqrt{n}}\sum_{i=1}^{N}
\left(\boldsymbol{H}^{\T}\left( \boldsymbol{\xi }_{n}\right) _{ji}\right)
\int_{-\infty }^{\infty }%
\left(
{\psi }\left( \boldsymbol{\eta}\left( \boldsymbol{\xi }_{n},z \right)\right) 
\right)_{i}
\phi\left( z\right) dz \\
&=&\frac{\kappa^2 }{n}\sum_{i=1}^{N}
\left(\boldsymbol{H}^{\T}\left( \boldsymbol{\xi }_{n}\right) _{ji}\right)
\int_{-\infty }^{\infty }%
\left(
{\psi }^{\prime}\left( \boldsymbol{\eta}\left( \boldsymbol{\xi }_{n},z \right)\right) 
\right)_{i}
\phi\left( z\right) dz\frac{\partial \boldsymbol{h}_{i}^{\T}\left( \boldsymbol{\xi }_{n}\right) 
\boldsymbol{\beta }}{\partial \beta_k} \\
&=&\frac{\kappa^2 }{n}\sum_{i=1}^{N}\left(\boldsymbol{H}^{\T}\left( \boldsymbol{\xi }_{n}\right) _{ji}\right)
\int_{-\infty }^{\infty }%
\left(
{\psi }^{\prime}\left( \boldsymbol{\eta}\left( \boldsymbol{\xi }_{n},z \right)\right) 
\right)_{i}
\phi\left( z\right) dzH_{i,k} \\
&=&\frac{\kappa^2 }{n}
\boldsymbol{H}^{\T} \left( \boldsymbol{\xi}_n \right)
\oplus _{i=1}^{N}  \Big[  \int_{-\infty }^{\infty }\psi^{\prime}\left(\boldsymbol{\eta}\left( \boldsymbol{\xi }_{n},z\right) \right)\phi \left( z\right) dz \Big]_{i}
\boldsymbol{H}\left( \boldsymbol{\xi}_n \right)
\end{eqnarray*}
as claimed. $\hfill \square $
\subsection{Proof of Lemma \protect\ref{lemma: maxdiff}}

We calculate that, using the Huber score function, 
\begin{equation*}
\frac{\sigma _{M}^{2}}{\sigma ^{2}}=\frac{1-2c\phi \left( c\right) +2\left(
c^{2}-1\right) \Phi \left( -c\right) }{\left( 1-2\Phi \left( -c\right)
\right) ^{2}}\overset{def}{=}G\left( c\right) ,
\end{equation*}%
in terms of which $\nu _{\text{\textsc{m}}}=\left( \gamma ^{2}G\left(
c\right) +1\right) ^{-1}$. The function $G\ $is the asymptotic variance of
the M-estimate of location with $N(0,1)$ errors. It is decreasing in $c$ and
must exceed the asymptotic variance of the sample mean: $G\left( c\right)
>G\left( \infty \right) =1$. Two applications of L'H\^{o}pital's rule give $%
G\left( 0\right) =\pi /2$, the asymptotic variance of the median of a normal
sample, with corresponding regression estimate given by $\psi \left(
x\right) =sign(x)$. Thus 
\begin{equation*}
0<\nu _{\text{\textsc{ls}}}-\nu _{\text{\textsc{m}}}=\frac{\gamma ^{2}}{%
\gamma ^{2}+1}\left\{ \frac{G\left( c\right) -1}{\gamma ^{2}G\left( c\right)
+1}\right\} \overset{def}{=}H\left( c,\gamma ^{2}\right) .
\end{equation*}%
For fixed $\gamma $, $H\left( c,\gamma ^{2}\right) $ is increasing in $%
G\left( c\right) $, so is maximized at $c=0$. Then $H\left( 0,\gamma
^{2}\right) $ vanishes at $\gamma ^{2}=0,\infty $ and is maximized at $%
\gamma ^{2}=\left[ G\left( 0\right) \right] ^{-1/2}$, with value $%
\max_{c,\gamma }\left( \nu _{\text{\textsc{ls}}}-\nu _{\text{\textsc{m}}%
}\right) =\left( \sqrt{G\left( 0\right) }-1\right) \left/ \left( \sqrt{%
G\left( 0\right) }+1\right) \right. $. At the maximum, $\nu _{\text{\textsc{%
ls}}}=1/\left( \gamma ^{2}+1\right) $ and $\nu _{\text{\textsc{m}}}=1/\left(
\gamma ^{2}G\left( 0\right) +1\right) $, giving (\ref{nus}) . \hfill\ $%
\square $

\medskip 

\section*{Acknowledgements}

This work was carried out with the support of the Natural Sciences and
Engineering Council of Canada. It has benefited from the insightful comments
of two anonymous referees.

\bibliographystyle{natbib}
\bibliography{references}

\end{document}